\documentclass[reqno,12pt]{amsart}
\usepackage[left=3cm, right=3cm, top=3cm, bottom=2cm]{geometry}
\usepackage{amsmath,amssymb}
\usepackage[mathscr]{euscript}
\usepackage{hyperref}
\usepackage{mathtools,multicol,enumitem,tikz,caption}

\newtheorem{theorem}{Theorem}[section]

\newtheorem{proposition}[theorem]{Proposition}
\newtheorem{corollary}[theorem]{Corollary}
\newtheorem{definition}[theorem]{Definition}
\newtheorem{rmk}[theorem]{Remark}
\numberwithin{equation}{section}

\begin{document}

\author[K. Kuoch]{Kevin Kuoch}  \address{Kevin Kuoch
  \hfill\break\indent MAP5, UMR CNRS 8145, Universit\'e Paris Descartes. 45 rue des Saints-P\`eres, 75006 Paris, France.}
\email{kevin.kuoch@gmail.com}

\title[Phase transition for a contact process with random slowdowns]{Phase transition for a contact process with random slowdowns}

\date{\today}

\begin{abstract}  
Motivated by a model of an area-wide integrated pest management, we develop an interacting particle system evolving in a random environment. It is a generalized contact process in which the birth rate takes two possible values, determined either by a dynamic or a static random environment. Our goal is to understand the phase diagram of both models by identifying the mechanisms that permit coexistence or extinction of the process.
\end{abstract}

\keywords{Interacting particle systems, contact process, percolation, phase transition, random environment.} 

\subjclass[2000]{Primary 60K35; Secondary 60K37, 92D25}

\maketitle
\thispagestyle{empty}

\section{Introduction}\label{sec:intro}
In this paper, we consider a generalized contact process describing the evolution on a lattice of a multitype population. Developed in the fifties by E. Knipling and R. Bushland (see \cite{KC05,Kni55}), the \textsl{Sterile insect technique} emerged to control the New World screw worm, a pest threatening livestocks in America. It is a pest control method whereby sterile individuals of the population, to either regulate or eradicate, are released. While sterile males compete with wild males, they eventually mate with (wild) females so that offsprings reach neither the adult phase nor sexual maturity, reducing the next generation. By repeated releases, we should be able to cause
a variety of outcomes ranging from reduction to extinction.

During the last decades, mosquito-borne diseases have become a serious worldwide spreading disease affecting millions of people each year. As the main vector of dengue fever, the mosquito-type \textsl{Aedes
aegypti}, adapted its lifestyle from the forests to urban areas across the world, it is considered as one of the most fastest growing disease. At short notice with no available cure, population control is the only viable option to stem the disease. 

Female mosquitoes taking a blood meal from a virus-infected person become carriers of the virus. After some days, the virus spreads to the mosquito's salivary glands and is therefore released into its saliva. From then on, this mosquito remains lifelong infected and transmits the virus to each new person it bites. As of today no cure exists, tackling the dengue transmission involves eradicating its vector: the mosquito. Lauded for several aspects and its success, the sterile insect technique is one of today main actors in this fight, as recent programs launched by Brazilian and European institutions highlight.

\bigskip
To understand this phenomenon, we construct a stochastic spatial model on a lattice in which the evolution of the wild population is governed by a contact process whose growth rate is slowed down in presence of sterile individuals, shaping a dynamic random environment.

The contact process with random slowdowns (CPRS) $(\eta_t)_{t\ge 0}$ has state space $\{0,1,2,3\}^{\mathbb Z^d}$. Each site of $\mathbb Z^d$ is either empty (state 0), occupied by wild individuals only (state 1), by sterile individuals only (state 2), or by both wild and sterile individuals (state 3). We only consider the type of the individuals present on each site (and not their number), it is biologically reasonable to assume enough females are around, not limiting matings. Sterile males are spontaneously released at rate $r$. The rate with which wild individuals give birth (to wild individuals) on neighbouring sites is $\lambda_1$ at sites in state 1 and $\lambda_2$ at sites in state 3. We assume that 
\begin{equation}\label{birth-rate}
\lambda_2<\lambda_1
\end{equation}
to highlight a competition between individuals occurs at sites in state 3 so that the fertility is decreased.  Deaths of each population occur at all sites at rate 1, they are mutually independent.  With no further restriction, this process is called \textsl{symmetric} while  in a so-called \textsl{asymmetric} case, sterile individuals stem births on sites they occupy.

A number of models with dynamic random environment has been investigated (see Broman \cite{Bro07}, Remenik \cite{Rem08}, Steif-Warfheimer \cite{SW08}, Garet-Marchand \cite{GM14}). On the other hand, identifying mechanisms yielding extinction or coexistence of species from ecology issues has been devoted in several papers (cf. Schinazi \cite{Sch97}, Durrett-Neuhauser \cite{DN97}) using multitype contact processes. The common feature is a type of individuals evolving as a contact process coexisting along with an other type evolving independently in a markovian way and dictating somehow its growth.

\bigskip
The question we address now is for which values of $r$ does the wild population survive or die out. Along with the monotonicity of the survival probability, we prove the existence and uniqueness of a phase transition with respect to the release rate $r$ for fixed growth rates $\lambda_1$ and $\lambda_2$. 

\bigskip
Denote by $\lambda_c(d)$ the critical value of the basic contact process on $\mathbb Z^d$. By comparing the process with basic contact processes, we explore in Propositions \ref{phase:trivial_sub} and \ref{phase:trivial_sup} varying growth parameters where no phase transition occurs i.e. when both $\lambda_1$ and $\lambda_2$ are smaller or larger than $\lambda_c(d)$. The most interesting cases are discussed in the following results, 
\begin{theorem}\label{transition:sym}
Suppose $\lambda_2 < \lambda_c(d) < \lambda_1$ fixed. Consider the symmetric CPRS. There exists a unique critical value $r_c \in (0,\infty)$ such that the wild population survives if $r < r_c$ and dies out if $r \geq r_c$.
\end{theorem}
\begin{theorem}\label{transition:asym}
Suppose $\lambda_c(d) < \lambda_1$ fixed. Consider the asymmetric CPRS. There exists a unique critical value $s_c \in (0,\infty)$ such that the wild population survives if $r < s_c$ and dies out if $r \geq s_c$.
\end{theorem}
These results agree with the former conclusions done by E. Knipling (1955) in a deterministic model. As a consequence, we discuss in Subsection \ref{transition:behaviour} the competitive ability of the sterile individuals: we prove the critical value increases as the competitiveness of the sterilized population decreases or as the fitness of the wild population increases.

To prove our results, we first construct the graphical representation of the process in Section \ref{sec:graph}, setting a percolation structure on which to define the process, lending itself to the use of percolation that we shall do in Section \ref{sec:phase}. Afterwards, we point out general properties of the system in Section \ref{sec:order}, such as necessary and sufficient conditions for the process to be monotone, then, only sufficient conditions in Propositions \ref{genprop:monotonicity_cns} and \ref{genprop:monotonicity_cs} to be in line with the construction of the process. It follows that the survival probability is monotone, see Corollary \ref{genprop:survivaldecrease}. The tricky part to prove these conditions lies in the definition of an order on the state space $\{0,1,2,3\}^{\mathbb Z^d}$, since an element of $\{0,1,2,3\}$ does not correspond to the number of particles but a type. 

Then, we use block constructions (see \cite{BD88}) and dynamic renormalization techniques (see e.g. \cite{BGN91}) to prove Theorems \ref{transition:sym} and \ref{transition:asym}  in Section \ref{sec:phase}.

\bigskip
So far, we were unable to get a hand on bounds for the critical rate. To this purpose, we first consider the mean-field equations in Section \ref{sec:mean}, featuring the densities of each type of individuals, we can explicitly find equilibria yielding numerical bounds on the phase transition. We derive a rigorous proof of the convergence of the empirical densities to these equations in \cite{kms}.

\bigskip
In a very different view, 
we consider the process $(\xi_t,\omega)_{t \geq 0}$, where $(\xi_t)_{t \geq 0}$ is a contact process on $\{0,1\}^{\mathbb Z}$ with varying growth rate determined by a quenched random environment $\omega$ initially fixed. In words, it describes an initial release of sterile individuals within the target area and looking at the induced time-evolution of the wild population. Using former results obtained by T.M. Liggett \cite{Lig91,Lig92}, we obtain in Section \ref{sec:quench} several survival and extinction conditions for the process. In that way, we consider two kinds of growth rates in $\mathbb Z$: one where the rates depend on the edges and one where the rates depend on the vertices. This yields numerical bounds on the phase transition for the process to survive or die out.

Random parameters for the contact process have become a matter of interest by showing up unexpected behaviours and interesting survival conditions of the process, see Andjel \cite{And92}, Bramson-Durrett-Schonmann \cite{BDS91}, Klein \cite{Kle94}, Newman-Volchan \cite{NV96}.

\bigskip
The paper is organized as follows. We begin by introducing the model and results for a phase transition in Section \ref{sec:settings}. In this view, we provide a graphical construction of the process in Section \ref{sec:graph} and stochastic order properties in Section \ref{sec:order}. Section \ref{sec:phase} is devoted to prove the phase transition on $\mathbb Z^d$. After what, we investigate the mean field-model in Section \ref{sec:mean} and the quenched model in Section \ref{sec:quench} both exhibiting bounds on the critical value.

\section{Settings and results}\label{sec:settings}

\subsection{The model}\label{subsec:model}
On the state space $\Omega = F^S$, where  $F=\{0,1,2,3\}$ and $S = \mathbb Z^d$, the \textsl{contact process with random slowdowns} (CPRS) is an interacting particle system $(\eta_t)_{ t \geq 0 }$ whose configuration at time $ t $ is $ \eta_t \in \Omega $, that is, for all $x \in \mathbb Z^d$, $\eta_t(x) \in F$ represents the state of site $x$ at time $t$. Two sites $x$ and $y$ are nearest neighbours on $\mathbb Z^d$ if $\|x-y \| = 1 $, also written $x \sim y $, and $n_i(x,\eta_t)$ stands for the number of nearest neighbours of $x$ in state $i$, $i=1,3$. We define now the symmetric and asymmetric CPRS we mentioned earlier as follows.

\begin{definition}[Asymmetric CPRS]
The asymmetric CPRS $(\eta_t)_{t \geq 0}$ has transition rates in $x$ for a current configuration $\eta$ :
\begin{equation}\label{asym-CPRS}
\begin{array}{l@{\extracolsep{2cm}}l}
0 \rightarrow 1 \ \mbox{ at rate } \lambda_1 n_1(x,\eta) + \lambda_2 n_3(x,\eta) & 1 \rightarrow 0 \ \mbox{ at rate } 1\\ 
0 \rightarrow 2 \ \mbox{ at rate } r & 2 \rightarrow 0 \ \mbox{ at rate  } 1\\ 
1 \rightarrow 3 \ \mbox{ at rate } r & 3 \rightarrow 1 \ \mbox{ at rate  } 1\\ 
& 3 \rightarrow 2 \ \mbox{ at rate } 1
\end{array}
\end{equation}
\end{definition}

\begin{definition}[Symmetric CPRS]
The symmetric CPRS $(\eta_t)_{t \geq 0}$ has transition rates in $x$ for a current configuration $\eta$ given by  \eqref{asym-CPRS}, to which one adds the following transition:
\begin{equation}
2 \rightarrow 3 \ \mbox{ at rate } \lambda_1 n_1(x,\eta) + \lambda_2 n_3(x,\eta).
\end{equation}
\end{definition}

Therefore, 
the 2's dictate independently the growth rate of the process, and even inhibit births in the asymmetric case, so that the 2's shape a \textsl{dynamic random environment} for the 1's.

\bigskip
In both cases, if $\eta \in \Omega$ and $x \in \mathbb Z^d$, denote by $\eta_x^i \in \Omega, \ \ i \in \{0,1,2,3\}$, the configuration obtained from $\eta$ after a flip of $ x $ to state $ i $:
\begin{equation}
 \eta \longrightarrow \eta_x^i \mbox{ at rate  $c(x,\eta,i)$,}
\mbox{ where } \forall u \in \mathbb Z^d,  \ \eta_x^i(u) = \left\{
\begin{array}{ll}
\eta(u) & \mbox{ if } u \neq x \\
i & \mbox{ if } u =x
\end{array}
\right.
\end{equation}
Let $\mathcal{L}$ be the infinitesimal generator of $ (\eta_t)_{t \geq 0} $, for any cylinder function $f$ on $\Omega$,
\begin{equation}\label{generateur1}
\mathcal{L}f(\eta) = \sum\limits_{x \in \mathbb Z^d} \sum\limits_{i=0}^3 c(x,\eta,i)\big(f(\eta_x^i)-f(\eta) \big)
\end{equation}
with infinitesimal transition rates, common to both cases (see \eqref{asym-CPRS}),
\begin{equation}\label{rates} \begin{array}{ll}
c(x,\eta,0) = 1 \mbox{ if } \eta(x) \in \{ 1, 2\},
& c(x,\eta,1) = \left\{
\begin{array}{l}
\lambda_1 n_1(x,\eta) + \lambda_2 n_3(x,\eta) \mbox{ if } \eta(x) = 0 \\
1 \mbox{ if } \eta(x) = 3 \\
\end{array}
\right., \\
c(x,\eta,2) = \left\{
\begin{array}{l}
r \mbox{ if } \eta(x) = 0 \\
1 \mbox{ if } \eta(x) = 3 \\
\end{array}
\right.,
& c(x,\eta,3) = 
\begin{array}{l}
r \mbox{ if } \eta(x) = 1 \\
\end{array} \\
\end{array} 
\end{equation}
and add the following rate in the symmetric case:
$$c(x,\eta,3) =
\lambda_1 n_1(x,\eta) + \lambda_2 n_3(x,\eta) \mbox{ if } \eta(x) = 2. $$
Since all the rates in \eqref{rates} are bounded, by \cite[Theorem 3.9]{Lig85} there exists a unique Markov process associated to the generator \eqref{generateur1}. Denote by $(\eta_t^A)_{t \geq 0}$ the process starting from $A$, i.e. such that $\eta_0 = \mathbf{1}_A $, in other words $\eta_0$ is the configuration with sites in state 1 in $A$ and empty otherwise.  We care about the evolution of the wild population, i.e. individuals contained in sites in state 1 and 3. Define
\begin{equation}\label{cpdre:occupied_set}
H_t^A =\{ x \in \mathbb Z^d : \eta_t^A(x) \in \{1 , 3\} \},
\end{equation}
as the set of sites containing the wild population at time $t \geq 0$. 

\bigskip
Denote by $\mathbb P_{\lambda_1, \lambda_2, r}$ the distribution of  $(\eta_t^{\{0\}})_{t \geq 0}$ with parameters  $(\lambda_1, \lambda_2, r)$.  When $\lambda_1$ and  $\lambda_2$ are fixed, we simply denote this law by $\mathbb P_r$ instead of $\mathbb P_{\lambda_1,\lambda_2,r}$.
The process $(\eta_t)_{t\geq 0}$ with initial configuration $ \eta_0 = \mathbf 1_{\{0\}}$ is said to \textsl{survive} if
\begin{equation}\label{cpsp:surviv.extinct}
\mathbb{P}_{\lambda_1, \lambda_2, r}(\forall t \geq 0, \ H_t^{\{0\}} \neq \emptyset) > 0
\end{equation}
and \textsl{die out} otherwise.

\bigskip
Remark when $r=0$, the process $(\eta_t)_{t \geq 0}$ is a basic contact process, denoted by $(\xi_t)_{t \geq 0}$, whose generator, deduced from \eqref{generateur1} with $r=0$, acts on local functions of $\Omega$ so that
$$\mathcal L^{cp}f(\xi) = \sum\limits_{x \in \mathbb Z^d} \Big\{ f(\xi_x^0)-f(\xi) + \lambda_1 \big(f(\xi_x^1) - f(\xi)\big) \Big\}.$$
It is known the basic contact process exhibits a unique phase transition on $\mathbb Z^d$, whose critical value is denoted by $\lambda_c(d)$, see \cite{Har74,Lig99,Lig85} for further details.

Define the \textsl{critical value} according to the parameter $r$ by
\begin{equation}\label{genprop:criticalval}
r_c = r_c (\lambda_1,\lambda_2) := \inf \{ r > 0 : \mathbb{P}_r(\exists t \geq 0, \ H_t^{\{0\}} = \emptyset) = 1 \}
\end{equation}
Indeed, the class $\{0,2\}$ is a trap: as soon as $H_t = \emptyset$, the wild population is extinct while sterile individuals (2) are infinitely often dropped along the time.

\subsection{Results}
As announced in Section \ref{sec:intro}, we begin by a first set of conditions for the process to survive or die out, when $\lambda_2 < \lambda_1$ are both smaller or larger than $\lambda_c$:
\begin{proposition}\label{phase:trivial_sub}
Suppose $\lambda_2 < \lambda_1 \leq \lambda_c(d)$. For all $r \geq 0$, both symmetric and asymmetric CPRS with parameters   $(\lambda_1,\lambda_2,r)$ die out.
\end{proposition}
\begin{proposition}\label{phase:trivial_sup}
Suppose $\lambda_c(d) < \lambda_2 < \lambda_1$. For all $r \geq 0$, the symmetric CPRS with parameters $(\lambda_1,\lambda_2,r)$ survives.
\end{proposition}

We now turn to Theorems \ref{transition:sym} and \ref{transition:asym}, for which we shall prove:
\begin{theorem}\label{transition:existence_sym}
Assume $\lambda_2 < \lambda_c < \lambda_1$ fixed. Let $(\eta_t)_{t \geq 0}$ be the symmetric CPRS. Then,
\begin{enumerate}[label={(\roman*)}]
\item there exists $r_0 \in (0,\infty)$ such that if $r < r_0$ then the process $(\eta_t)_{t \geq 0}$ survives.
\item there exists $r_1 \in (0,\infty)$ such that if $r > r_1$ then the process $(\eta_t)_{t \geq 0}$ dies out.
\end{enumerate}
\end{theorem}
\begin{theorem}\label{transition:existence_asym}
Assume $\lambda_c(d) < \lambda_1$ and $\lambda_2 < \lambda_1$ fixed. Let $(\eta_t)_{t \geq 0}$ be the asymmetric CPRS. Then,
\begin{enumerate}[label={(\roman*)}]
\item there exists $s_0 \in (0,\infty)$ such that if $r < s_0$ then the process $(\eta_t)_{t \geq 0}$ survives.
\item there exists $s_1 \in (0,\infty)$ such that if $r > s_1$ then the process $(\eta_t)_{t \geq 0}$ dies out.
\end{enumerate}
\end{theorem}

Following Theorems \ref{transition:existence_sym} and \ref{transition:existence_asym}, monotonicity arguments (given by Corollary \ref{genprop:survivaldecrease}) yield the existence of the critical value $r_c$ \eqref{genprop:criticalval}. We complete the study by investigating the behaviour of the critical process.

\begin{theorem}\label{transition:diesout} In both cases,
the critical CPRS dies out: $\mathbb P_{r_c} \big( \forall t \geq 0, \ H_t^{\{0\}} \neq \emptyset \big) = 0.$
\end{theorem}
We shall prove these results in Section \ref{sec:phase}.

\section{Graphical representation}\label{sec:graph}
In view of using percolation, we construct the process from a collection of independent Poisson processes, see \cite{Har78}. Think of the diagram $\mathbb Z^d \times \mathbb R_+$. For each $x \in \mathbb Z^d$, consider the arrival times of mutually independent families of Poisson processes: $\{ A^{x}_n: n \geq 1 \}$ with rate $r$, $\{ D^{1,x}_n: n \geq 1 \}$ and $\{ D^{2,x}_n: n \geq 1 \}$ with rate $1$ and for any $y$ such that $y \sim x$, $\{ T^{ x,y}_n: n \geq 1 \}$ with rate $\lambda_1$. Let  $\{U_n^{x}: n \geq 1\}$ be independent uniform random variables on $(0,1)$, independent of the Poisson processes.

\bigskip
At space-time point $(x, A^x_n)$, put a ``\begin{tikzpicture}
\filldraw[draw=black] circle (2pt);
\end{tikzpicture}" to indicate, if $x$ is occupied by type-1 individuals (resp. empty), that it turns into state 3 (resp. state 2). 
At $(x, D^{1,x}_n)$ (resp. at $(x, D^{2,x}_n)$), put an ``X" (resp.``\begin{tikzpicture}
\filldraw[fill=white,draw=black] circle (2pt);
\end{tikzpicture}") to indicate at $x$, that a death of type-1
(resp. of type-2)
occurs. At times $T^{x,y}_n$, draw an arrow from $x$ to $y$ and two kinds of actions occur according to the occupation at $x$:
if $x$ is occupied by type-1 individuals, the arrow indicates a birth in $y$ of a type-1 individual if $y$ is empty or in state 2 
; if $x$ is occupied by type-3 individuals giving birth at rate $\lambda_2 < \lambda_1$, check at $(x,T_n^{x,y})$ if $U_n^{x} < \lambda_2 / \lambda_1$ to indicate, if success, that the arrow is effective so that a birth in $y$ of a type-1 individual occurs if $y$ is empty or in state 2
. In the asymmetric case, births occur only if $y$ is empty.

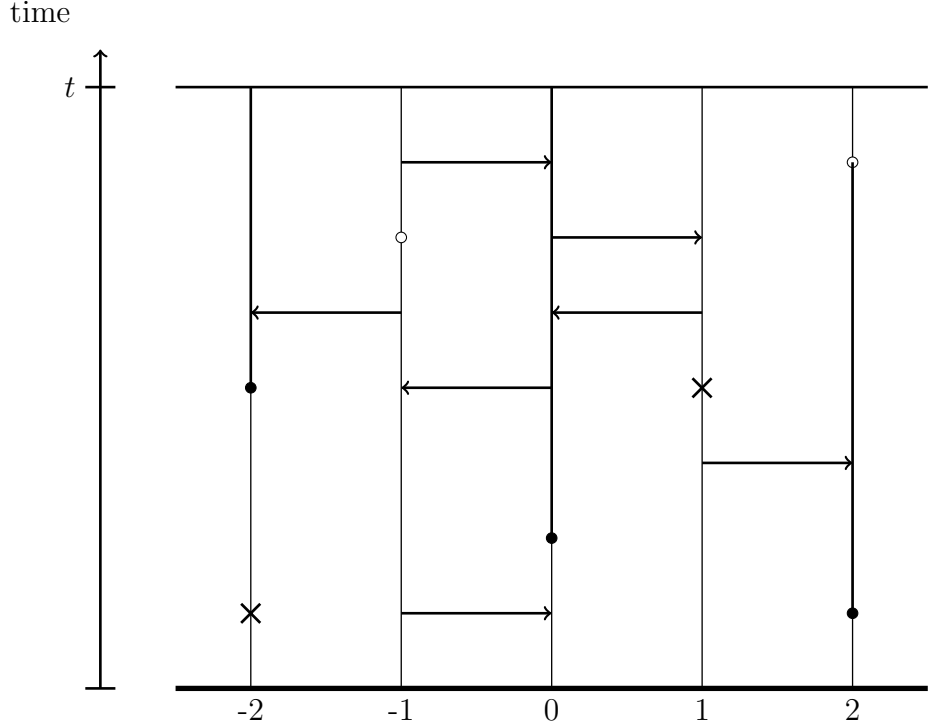
\begin{figure}[h]
\begin{center}
\begin{tikzpicture}
\draw[help lines,line width=2pt,step=2,draw=black] (-4,0) to (6,0);
\draw[help lines,line width=1pt,step=2,draw=black] (-4,8) to (6,8);
\draw[help lines,line width=0.5pt,step=2,draw=black] (-3,0) to (-3,8);
\draw[help lines,line width=0.5pt,step=2,draw=black] (-1,0) to (-1,8);
\draw[help lines,line width=0.5pt,step=2,draw=black] (1,0) to (1,8);
\draw[help lines,line width=0.5pt,step=2,draw=black] (3,0) to (3,8);
\draw[help lines,line width=0.5pt,step=2,draw=black] (5,0) to (5,8);

\node[anchor=north] at (-3,0) {-2};
\node[anchor=north] at (-1,0) {-1};
\node[anchor=north] at (1,0) {0};
\node[anchor=north] at (3,0) {1};
\node[anchor=north] at (5,0) {2};

\draw[->,line width=1pt,black] (-1,1)--(1,1)  node [midway,above] {};
\draw[->,line width=1pt,black] (3,3)--(5,3)  node [midway,above] {};
\draw[->,line width=1pt,black] (-1,5)--(-3,5)  node [midway,above] {};
\draw[->,line width=1pt,black] (3,5)--(1,5)  node [midway,above] {};
\draw[->,line width=1pt,black] (1,4)--(-1,4)  node [midway,above] {};
\draw[->,line width=1pt,black] (-1,7)--(1,7)  node [midway,above] {};
\draw[->,line width=1pt,black] (1,6)--(3,6)  node [midway,above] {};

\filldraw[fill=white,draw=black] (-1,6) circle [radius=2pt];
\node at (3,4) {$\mathbin{\tikz [x=1.4ex,y=1.4ex,line width=.2ex,black] \draw (0,0) -- (1,1) (0,1) -- (1,0);}$};
\filldraw[fill=white,draw=black] (5,7) circle [radius=2pt];
\node at (-3,1) {$\mathbin{\tikz [x=1.4ex,y=1.4ex,line width=.2ex,black] \draw (0,0) -- (1,1) (0,1) -- (1,0);}$};
\draw[fill,black] (1,2) circle [radius=2pt];
\draw[fill,black] (5,1) circle [radius=2pt];
\draw[fill,black] (-3,4) circle [radius=2pt];
\draw[help lines,line width=1pt,step=2,draw=black] (-3,4) to (-3,8);
\draw[help lines,line width=1pt,step=2,draw=black] (1,2) to (1,8);
\draw[help lines,line width=1pt,step=2,draw=black] (5,1) to (5,7);


\draw [->,line width=1pt] (-5,0) -- (-5,8.5);
\node at (-5.8,9) {$\mbox{time}$};
\draw[help lines,line width=1pt,step=2,draw=black] (-5.2,0) to (-4.8,0);
\draw[help lines,line width=1pt,step=2,draw=black] (-5.2,8) to (-4.8,8);
\node at (-5.4,8) {$t$};
\end{tikzpicture}
\end{center}
\caption{Starting from $\eta_0 = \mathbf{1}_{\{0\}}$, following the arrows, if $U_1^0 < \frac{\lambda_2}{\lambda_1}$ and $U_2^0 < \frac{\lambda_2}{\lambda_1}$, the wild population occupies at time $t$ the set $H_t = \{-1,0,1\}$ in the asymmetric case and the set $H_t=\{-2,-1,0,1\}$ in the symmetric case.}
\label{fig:graphrep}
\end{figure}

See Figure \ref{fig:graphrep} for an example of the time-evolution of both processes starting from an identical initial configuration, in the space-time picture $\mathbb Z \times \mathbb R_+$.\\

For $s \leq t$, there exists an \textsl{active path from $(x,s)$ to $(y,t)$} in $\mathbb Z^d \times \mathbb R_+$ is there exists a sequence of times $s=s_0 < s_1 < ... < s_{n-1} < s_n = t$ and a sequence of corresponding spatial locations $x=x_0, x_1, ..., x_n = y$ such that:
\begin{enumerate}[label=\roman{*}.]
\item for $i = 1, ..., n-1 $, vertical segments $\{ x_i \} \times (s_i, s_{i+1})$ do not contain any X's.
\item for $i = 1, ..., n $, there is an arrow from $x_{i-1}$ to $x_i$ at times $s_i$ and if $x_{i-1} \times s_i$ is lastly preceded by a `` \begin{tikzpicture}
\filldraw[draw=black] circle (2pt);
\end{tikzpicture}"
this arrow exists only if $U_{s_i}^{x_{i-1}} < \lambda_2 / \lambda_1$.
\end{enumerate}

and in the asymmetric case, substitute ii. by
\begin{enumerate}[label=ii'.]
\item for $i = 1, ..., n $, there is an arrow from $x_{i-1}$ to $x_i$ at times $s_i$ while $\{x_i\}\times s_i$ is not lastly preceded by a ``\begin{tikzpicture}
\filldraw[draw=black] circle (2pt);
\end{tikzpicture}'', while if $x_{i-1} \times s_i$ is lastly preceded by a  ``\begin{tikzpicture}
\filldraw[draw=black] circle (2pt);
\end{tikzpicture}''
this arrow is effective if $U_{s_i}^{x_{i-1}} < \lambda_2 / \lambda_1$.
\end{enumerate}

\bigskip
Consider the process $(A_t^A)_{t \geq 0}$, the set of sites at time $t$ reached by active paths starting from an initial configuration $A_0=A$, containing sites in state 1 in $A$ and 0 otherwise:
$$A_t^A = \{ y \in \mathbb Z : \exists x \in A \mbox{ such that } (x,0) \rightarrow (y,t) \}$$
Then $A_t^A = H_t^A$, with $H_t^A$ defined in \eqref{cpdre:occupied_set} so that $A_t^A$ represents the wild population at time $t$ starting from an initial configuration $A$ of type-1 individuals.

\bigskip
This graphical representation allows us to couple CPRS starting from different initial configurations by imposing the evolution to obey to the same Poisson processes. Other kinds of couplings are possible through the analytical construction of the process as we will see. 
More generally, graphical representations allow to couple processes with different dynamics as well, we investigate this in the next section.

\section{Stochastic order and attractiveness}\label{sec:order}
The stochastic order between processes is related to the total order defined on the set of values $F = \{ 0, 1, 2, 3\}$. In a biological context, setting an order between types of individuals does not make sense, but mathematically it allows us to construct a monotone model and to compare different dynamics. 
To avoid confusion, denote by $A$ the state 2, by $B$ the state 0, by $C$ the state 3 and by $D$ the state 1. For two given configurations $\eta^{(2)}, \eta^{(1)} \in \Omega$, we have $\eta^{(1)} \leq \eta^{(2)}$, if coordinate-wise
$$\eta^{(1)}(x) \leq \eta^{(2)}(x) \ \mbox{ for all } x \in \mathbb Z^d, $$
with respect to the following order on $F$:
\begin{equation}\label{ordreF}
A < B < C < D
\end{equation}
It is important to note that this order is the only one possible with which such a stochastic order for the the process $(\eta_t)_{t \geq 0}$ holds.

We now settle necessary and sufficient conditions, then only sufficient, to obtain properties of stochastic order which give monotonicity properties and comparisons between processes. 
Using Borrello's results \cite{Bor11}, let $x, y \in \mathbb Z^d$ be two neighbouring sites and $\alpha, \beta \in F$ to rewrite the transition rates for $k \in \{1,2\}$ thanks to
\begin{enumerate}[label=$\bullet$]
\item $R_{\alpha, \beta}^{0,k} $ is the birth rate of a type-1 individual in $y$ when $\eta(y) = \beta$ and depending only on the value of $\eta$ in $x$ given by $\alpha$. The state in $y$ flips from $\beta$ to $\beta+k$.
\item $P_{\beta}^{k}$ is the jump rate of a site in state $\beta$ to state $\beta + k$.
\item $P_{\alpha}^{- k}$ is the jump rate of a site in state $\alpha$ to state $\alpha - k$ for $k \leq \alpha$.
\end{enumerate}
 \begin{theorem}\label{ips:borrello}\cite[Theorem 2.4]{Bor11}
For all $(\alpha,\beta) \in F^2$, $(\gamma,\delta) \in F^2 $ such that  $(\alpha, \beta) \leq (\gamma,\delta)$ (coordinate-wise i.e. $\alpha \leq \gamma$ and $\beta \leq \delta$), $h_1 \geq 0$, $j_1 \geq 0$, an interacting particle systems $(A_t)_{t \geq 0}$ with transition rates $(R_{\alpha, \beta}^{0,k}, P_{\beta}^{+ k}, P_{\alpha}^{-k})$ is stochastically larger than an interacting particle system $(B_t)_{t \geq 0 }$ with transition rates $(\widetilde{R}_{\alpha, \beta}^{0,k}, \widetilde{P}_{\beta}^{+ k}, \widetilde P_{\alpha}^{-k})$ if and only if
\begin{equation}\label{borrello}i) \sum\limits_{k>\delta-\beta+j_1} \widetilde{\Pi}_{\alpha,\beta}^{0,k} \leq \sum\limits_{l > j_1} \Pi_{\gamma,\delta}^{0,l} \ \mbox{ and } \ ii) \sum\limits_{k> h_1} \widetilde{\Pi}_{\alpha,\beta}^{-k,0} \geq \sum\limits_{l > \gamma - \alpha +h_1} \Pi_{\gamma,\delta}^{-l,0} 
\end{equation}
where $\Pi^{0,k}_{\alpha,\beta} = R^{0,k}_{\alpha,\beta} + P^k_\beta \mbox{ and } \Pi^{-k,0}_{\alpha,\beta} = P^{-k}_\alpha$.
\end{theorem}
One has for the asymmetric CPRS, the following rates.
\begin{equation}\label{borrello_CPDRE}
\begin{array}{cc}
R_{D,B}^{0,2} = \lambda_1, \ \ R_{C,B}^{0,2} = \lambda_2, \ \ P_A^1 = P_C^1 = 1, \\
P_B^{-1} =P_D^{-1} = r, \ \ P_C^{-2} = P_D^{-2} = 1,
\end{array}
\end{equation}
to which, one adds the following rates if we consider the symmetric CPRS,
\begin{equation}\label{borrello_CPDRE2}
R_{D,A}^{0,2} = \lambda_1, \ \ R_{C,A}^{0,2} = \lambda_2.
\end{equation}
Similarly, for a basic contact process with growth rate $\lambda_1$ on $\{0,1\}^{\mathbb Z^d}$, one has
\begin{equation}\label{borrello_contactsuper}
\widetilde{R}_{D,B}^{0,2} = \lambda_1, \ \widetilde{P}_{D}^{-2} = 1.
\end{equation}
For a basic contact process $(\widetilde \xi_t)_{t \geq 0}$ with growth rate $\lambda_2$, defined on $\{2,3\}^{\mathbb Z^d}$, let
\begin{equation}\label{borrello_contactsub}
\widetilde{R}_{C,A}^{0,2} = \lambda_2, \ \widetilde{P}_{C}^{-2} = 1.
\end{equation}
Anticipating forthcoming computations, we expand \eqref{borrello} when $(A_t)_{t \geq 0}$ and $(B_t)_{t \geq 0}$ correspond to two symmetric CPRS with respective rates $(\lambda_1,\lambda_2,r)$ denoted with upper indexes $(2)$ and $(1)$. Note that by \eqref{borrello_CPDRE}-\eqref{borrello_contactsub}, the following inequalities can be used as well for asymmetric CPRS or basic contact processes. Taking $j_1, h_1 \in \{0,1\}$ (one can check they are the only non trivial possible values), we have respectively from \eqref{borrello}(i)-(ii) applied to our model:

\begin{footnotesize}
\begin{align}\label{genprop:monotonicity_cns1} \nonumber
& \begin{multlined} \mathbf{1}\{j_1=0\} \mathbf{1}\{k=2\} \mathbf{1}\{\delta-\beta=1\}  \Big( \mathbf{1}\{\delta=C,\beta=B\} \big( R_{D,B}^{0,2,(1)}\mathbf{1}\{\alpha=D\} + R_{C,B}^{0,2,(1)}\mathbf{1}\{\alpha=C\} \big) \\ 
+ \mathbf{1}\{\delta=B,\beta=A\} \big( R_{D,A}^{0,2,(1)} \mathbf{1}\{\alpha=D\} + R_{C,A}^{0,2,(1)} \mathbf{1}\{\alpha=C\}\big) \Big) \\
+ \mathbf{1}\{j_1=0\} \mathbf{1}\{k=2\} \mathbf{1}\{\delta-\beta=0\}  \Big( \mathbf{1}\{\delta=\beta=B\} \big( R_{D,B}^{0,2,(1)}\mathbf{1}\{\alpha=D\} + R_{C,B}^{0,2,(1)}\mathbf{1}\{\alpha=C\} \big) \\ 
+ \mathbf{1}\{\delta=\beta=A\} \big( R_{D,A}^{0,2,(1)} \mathbf{1}\{\alpha=D\} + R_{C,A}^{0,2,(1)} \mathbf{1}\{\alpha=C\} \big) \Big)   \\
+ \mathbf{1}\{j_1=0\} \mathbf{1}\{k=1\} \mathbf{1}\{\delta-\beta=0\} \Big( \mathbf{1}\{\delta=\beta=C\} P_C^{1,(1)} +\mathbf{1}\{\delta=\beta=A\} P_A^{1,(1)} \Big) \\
+ \mathbf{1}\{j_1=1\} \mathbf{1}\{k=2\} \mathbf{1}\{\delta-\beta=0\}  \\
\Big( \mathbf{1}\{\delta=\beta=B\} \big( R_{D,B}^{0,2,(1)}\mathbf{1}\{\alpha=D\} + R_{C,B}^{0,2,(1)}\mathbf{1}\{\alpha=C\} \big) \\
+ \mathbf{1}\{\delta=\beta=A\} \big( R_{D,A}^{0,2,(1)}\mathbf{1}\{\alpha=D\} + R_{C,A}^{0,2,(1)} \mathbf{1}\{\alpha=C\}\big) \Big) \end{multlined} \\
& \begin{multlined} \leq  \mathbf{1}\{j_1=0\} \mathbf{1}\{l=2\}  \Big( \mathbf{1}\{\delta=B\} \big( R_{D,B}^{0,2,(2)}\mathbf{1}\{\gamma=D\} + R_{C,B}^{0,2,(2)}\mathbf{1}\{\gamma=C\} \big) \\ 
+ \mathbf{1}\{\delta=A\} \big( R_{D,A}^{0,2,(2)}\mathbf{1}\{\gamma=D\} + R_{C,A}^{0,2,(2)} \mathbf{1}\{\gamma=C\}\big) \Big)\\
 + \mathbf{1}\{j_1=0\} \mathbf{1}\{l=1\} \Big( \mathbf{1}\{\delta=C\} P_C^{1,(2)}+ \mathbf{1}\{\delta=A\} P_A^{1,(2)} \Big) \\
+  \mathbf{1}\{j_1=1\} \mathbf{1}\{l=2\} \Big( \mathbf{1}\{\delta=B\} \big( R_{D,B}^{0,2,(2)} \mathbf{1}\{\gamma=D\}+ R_{C,B}^{0,2,(2)} \mathbf{1}\{\gamma=C\}\big) \\ + \mathbf{1}\{\delta=A\} \big( R_{D,A}^{0,2,(2)} \mathbf{1}\{\gamma=D\}+ R_{C,A}^{0,2,(2)} \mathbf{1}\{\gamma=C\}\big) \Big) \end{multlined} 
\end{align}
\end{footnotesize}
and
\begin{footnotesize}
\begin{align}\label{genprop:monotonicity_cns2} \nonumber
& \begin{multlined} \mathbf{1}\{h_1=0\} \mathbf{1}\{k=1\} \Big( \mathbf{1}\{\alpha=D\} P_D^{-1,(1)} + \mathbf{1}\{\alpha=B\} P_B^{-1,(1)} \Big) \\
+\mathbf{1}\{h_1=0\} \mathbf{1}\{k=2\} \Big( \mathbf{1}\{\alpha=D\} P_D^{-2,(1)} + \mathbf{1}\{\alpha=C\} P_C^{-2,(1)} \Big) \\
+ \mathbf{1}\{h_1=1\} \mathbf{1}\{k=2\} \Big( \mathbf{1}\{\alpha=D\} P_D^{-2,(1)} + \mathbf{1}\{\alpha=C\} P_C^{-2,(1)} \Big) \end{multlined} \\
& \begin{multlined}\geq \mathbf{1}\{h_1=0\} \mathbf{1}\{l=2\} \mathbf{1}\{\gamma-\alpha=1\} \Big( \mathbf{1}\{\gamma=D,\alpha=C\} P_D^{-2,(2)} + \mathbf{1}\{\gamma=C,\alpha=B\} P_C^{-2,(2)} \big) \\ + \mathbf{1}\{h_1=0\} \mathbf{1}\{l=2\} \mathbf{1}\{\gamma-\alpha=0\} \Big( \mathbf{1}\{\gamma=\alpha=D\}P_D^{-2,(2)} + \mathbf{1}\{\gamma=\alpha=C\} P_C^{-2,(2)}\Big) 
\\ + \mathbf{1}\{h_1=0\} \mathbf{1}\{l=1\} \mathbf{1}\{\gamma-\alpha=0\} \Big( \mathbf{1}\{\gamma=\alpha=D\}P_D^{-1,(2)} + \mathbf{1}\{\gamma=\alpha=C\} P_C^{-1,(2)}\Big) \\
+\mathbf{1}\{h_1=1\} \mathbf{1}\{l=2\} \mathbf{1}\{\gamma-\alpha=0\} 
\Big( \mathbf{1}\{\gamma=\alpha=D\} P_D^{-2,(2)} + \mathbf{1}\{\alpha=\gamma=C\} P_C^{-2,(2)} \Big)
\end{multlined}.
\end{align}
\end{footnotesize}
We are now ready to state and prove the following results.
\begin{proposition}\label{genprop:monotonicity_cns}
The symmetric CPRS is monotone, in the sense that one can construct on a same probability space two symmetric CPRS $(\eta_t^{(1)})_{t \geq 0}$ and $(\eta_t^{(2)})_{t \geq 0}$ with respective parameters  $(\lambda_1^{(1)}, \lambda_2^{(1)},r^{(1)})$ and $(\lambda_1^{(2)}, \lambda_2^{(2)},r^{(2)})$, such that
\begin{equation}
\eta^{(1)}_0 \leq \eta^{(2)}_0 \Longrightarrow \eta^{(1)}_t \leq \eta^{(2)}_t \mbox{ a.s. for all } t \geq 0
\end{equation}
if and only if all parameters satisfy 
\begin{multicols}{4}
\begin{enumerate}[label={\arabic*.}]
\item $\lambda_2^{(1)} \leq \lambda_1^{(1)}$, 
\item $\lambda_2^{(2)} \leq \lambda_1^{(2)}$,
\item $\lambda_1^{(1)} \leq \lambda_1^{(2)}$, 
\item $\lambda_2^{(1)} \leq \lambda_2^{(2)}$,
\item $r^{(1)} \geq r^{(2)}$
\item $\lambda_1^{(1)} \leq 1$, 
\item $\lambda_2^{(1)}  \leq 1$, 
\item $r^{(1)} \geq 1.$
\end{enumerate}
\end{multicols}
\end{proposition}
Remark conditions 1. and 2. are given by assumption \eqref{birth-rate}. 
\begin{proof}
Let $(\eta_t^{(1)})_{t \geq 0}$ and $(\eta_t^{(2)})_{t \geq 0}$ be two symmetric CPRS with respective parameters $(\lambda_1^{(1)}, \lambda_2^{(1)},r^{(1)})$ and $(\lambda_1^{(2)}, \lambda_2^{(2)},r^{(2)})$. 
By Theorem \ref{ips:borrello}, 
necessary and sufficient conditions on the rates for $(\eta_t^{(2)})_{t \geq 0}$ to be stochastically larger than $(\eta_t^{(1)})_{t \geq 0}$ are given 
with $(\alpha,\beta) \leq (\gamma,\delta)$, by \eqref{genprop:monotonicity_cns1}-\eqref{genprop:monotonicity_cns2}.
All different possible scenarios provide the following necessary conditions:
\begin{enumerate}[label={(\Roman*)}]
\item $j_1 \in \{0,1\}$, $\delta=\beta \in \{A,B\}$ in \eqref{genprop:monotonicity_cns1} give
\begin{enumerate}[label={(\roman*)}]
\item $\alpha=C, \gamma=D$: $\lambda_2^{(1)} \leq \lambda_1^{(2)}$. This is a consequence of conditions 1. and 3. or 2. and 4.
\item $\alpha=\gamma=C$: $\lambda_2^{(1)} \leq \lambda_2^{(2)}$ stated by condition 4.
\item $\alpha=\gamma=D$: $\lambda_1^{(1)} \leq \lambda_1^{(2)}$ stated by  condition 3.
\end{enumerate}
\item $j_1=0$, $\beta=B$, $\delta=1+\beta=C$ in \eqref{genprop:monotonicity_cns1} give
\begin{enumerate}[label={(\roman*)}]
\item $\alpha=D$: $\lambda_1^{(1)} \leq 1 $ stated by condition 6.
\item$\alpha=C$: $\lambda_2^{(1)} \leq 1 $ stated by condition 7.
\end{enumerate}
\item $h_1=0$, $\gamma=\alpha \in \{B,D\}$ in \eqref{genprop:monotonicity_cns2} give $r^{(1)} \geq r^{(2)}$ stated by condition 5.
\item $h_1=0$,$\alpha=B,\gamma=1+\alpha=C$ in \eqref{genprop:monotonicity_cns2} give $r^{(1)} \geq  1$ stated by condition 8.
\end{enumerate}
while in other scenarios, one retrieves redundantly the above conditions or tautological inequalities such as ``$1 \geq 0$". Finally, we obtained the conditions stated 1. to 8.
\end{proof}

Following Proposition \ref{genprop:monotonicity_cns}, we now construct a coupled process $(\eta_t^{(1)},\eta_t^{(2)})_{t \geq 0}$ on $\Omega \times \Omega$ such that $\eta_0^{(1)} \leq \eta_0^{(2)}$. According to the given order \eqref{ordreF} on $F$, as  $\eta_0^{(1)} \leq \eta_0^{(2)}$: 
$$ n_1(x,\eta_0^{(1)}) + n_3(x,\eta_0^{(1)}) \leq n_1(x,\eta_0^{(2)}) + n_3(x,\eta_0^{(2)}).$$ 
In what follows, we construct the coupling through generators. The three following tables depict the infinitesimal transitions of the coupled process.
\begin{center}
$
\begin{footnotesize}
\begin{array}{|c|c|}
\hline
\mbox{transition} & \mbox {rate}\\
   \hline
(0,0) \longrightarrow \left\{ \begin{array}{ll}
(1,1) \\(0,1) \\ (2,2) \\ (2,0)
\end{array}\right.
& \begin{array}{cc}
\lambda_1^{(1)} n_1(x,\eta^{(1)}) + \lambda_2^{(1)} n_3(x,\eta^{(1)}) \\
\lambda_1^{(2)}n_1(x,\eta^{(2)})- \lambda_1^{(1)} n_1(x,\eta^{(1)}) + \lambda_2^{(2)} n_3(x,\eta^{(2)})- \lambda_2^{(1)}n_3(x,\eta^{(1)})\\r^{(2)} \\r^{(1)} - r^{(2)}
\end{array}\\
(1,1) \longrightarrow \left\{ \begin{array}{ll}
(0,0) \\ (3,3) \\ (3,1)
\end{array}\right.
& \begin{array}{cc}
1 \\r^{(2)} \\ r^{(1)} - r^{(2)} 
\end{array}\\
(2,2) \longrightarrow \left\{ \begin{array}{ll}
(3,3) \\(2,3) \\ (0,0)
\end{array}\right.
& \begin{array}{cc}
\lambda_1^{(1)} n_1(x,\eta^{(1)}) + \lambda_2^{(1)} n_3(x,\eta^{(1)}) \\
\lambda_1^{(2)}n_1(x,\eta^{(2)})- \lambda_1^{(1)} n_1(x,\eta^{(1)}) + \lambda_2^{(2)} n_3(x,\eta^{(2)})- \lambda_2^{(1)}n_3(x,\eta^{(1)})\\1
\end{array}\\
(3,3) \longrightarrow \left\{ \begin{array}{ll}
(1,1) \\ (2,2)
\end{array}\right.
& \begin{array}{cc}
1 \\ 1
\end{array}\\
(2,0) \longrightarrow \left\{ \begin{array}{ll}
(3,1)\\(2,1) \\ (0,0) \\ (2,2)
\end{array}\right.
& \begin{array}{cc}
\lambda_1^{(1)} n_1(x,\eta^{(1)}) + \lambda_2^{(1)} n_3(x,\eta^{(1)}) \\
\lambda_1^{(2)}n_1(x,\eta^{(2)})- \lambda_1^{(1)} n_1(x,\eta^{(1)}) + \lambda_2^{(2)} n_3(x,\eta^{(2)})- \lambda_2^{(1)}n_3(x,\eta^{(1)})\\ 1 \\ r^{(2)} 
\end{array}\\
(2,3) \longrightarrow \left\{ \begin{array}{ll}
(3,3)\\(0,1) \\ (2,2)
\end{array}\right.
& \begin{array}{cc}
\lambda_1^{(1)} n_1(x,\eta^{(1)}) + \lambda_2^{(1)} n_3(x,\eta^{(1)}) \\1 \\ 1
\end{array}\\
(2,1) \longrightarrow \left\{ \begin{array}{ll}
(3,1)\\(2,0) \\ (0,1) \\ (2,3)
\end{array}\right.
& \begin{array}{cc}
\lambda_1^{(1)} n_1(x,\eta^{(1)}) + \lambda_2^{(1)} n_3(x,\eta^{(1)}) \\1\\1 \\ r^{(2)}
\end{array}\\
(3,1) \longrightarrow \left\{ \begin{array}{ll}
(2,0) \\ (1,1) \\ (3,3)
\end{array}\right.
& \begin{array}{cc}
1\\1 \\ r^{(2)}
\end{array}\\
(0,1) \longrightarrow \left\{ \begin{array}{ll}
(2,3) \\ (1,1) \\ (2,1) \\ (0,0)
\end{array}\right.
& \begin{array}{cc}
 r^{(2)} \\ \lambda_1^{(1)}n_1(x,\eta^{(1)}) + \lambda_2^{(1)} n_3(x,\eta^{(1)}) \\ r^{(1)}-r^{(2)} \\ 1
\end{array}\\
\hline
\end{array}
\end{footnotesize}
$
\captionof{table}{Coupled transitions of two symmetric CPRS, via a basic coupling, that preserve the partial order.}
\label{genprop:monotonicity_cns_tab1}
\end{center}

\begin{center}
$
\begin{footnotesize}
\begin{array}{|c|c|}
\hline
\mbox{transition} & \mbox {rate}\\
   \hline
(0,3) \longrightarrow \left\{ \begin{array}{ll}
(1,1)\\(0,1)\\(2,2)\\(2,3)
\end{array}\right.
& \begin{array}{cc}
\lambda_1^{(1)}n_1(x,\eta^{(1)}) + \lambda_2^{(1)} n_3(x,\eta^{(1)}) \\ (1-\lambda_1^{(1)})n_1(x,\eta^{(1)}) +(1- \lambda_2^{(1)}) n_3(x,\eta^{(1)}) \\ 1 \\ r^{(1)}-1
\end{array}\\
\hline
\end{array}
\end{footnotesize}
$
\captionof{table}{Coupled transitions of two symmetric CPRS, not via a basic coupling, that preserve the partial order.}
\label{genprop:monotonicity_cns_tab2}
\end{center}

\begin{center}
$
\begin{footnotesize}
\begin{array}{|c|c|}
\hline
\mbox{transition} & \mbox {rate}\\
   \hline
(1,0) \longrightarrow \left\{ \begin{array}{ll}
(1,1) \\ (3,2) \\ (3,0) \\ (0,0)
\end{array}\right.
& \begin{array}{cc}
\lambda_1^{(2)}n_1(x,\eta^{(2)}) + \lambda_2^{(2)} n_3(x,\eta^{(2)}) \\r^{(2)}\\ r^{(1)} - r^{(2)} \\ 1
\end{array}\\
(0,2) \longrightarrow \left\{ \begin{array}{ll}
(1,3) \\ (0,3) \\ (0,0) \\ (2,2)
\end{array}\right.
& \begin{array}{cc}
\lambda_1^{(1)}n_1(x,\eta^{(1)}) + \lambda_2^{(1)} n_3(x,\eta^{(1)}) \\\lambda_1^{(2)}n_1(x,\eta^{(2)})- \lambda_1^{(1)} n_1(x,\eta^{(1)}) + \lambda_2^{(2)} n_3(x,\eta^{(2)})- \lambda_2^{(1)}n_3(x,\eta^{(1)}) \\ 1 \\ r^{(1)}
\end{array}\\
(1,2) \longrightarrow \left\{ \begin{array}{ll}
(1,3)\\(3,2) \\ (1,0) \\ (0,2)
\end{array}\right.
& \begin{array}{cc}
\lambda_1^{(2)}n_1(x,\eta^{(2)}) + \lambda_2^{(2)} n_3(x,\eta^{(2)}) \\r^{(1)} \\ 1 \\ 1
\end{array}\\
(1,3) \longrightarrow \left\{ \begin{array}{ll}
(0,2) \\ (3,3) \\ (1,1)
\end{array}\right.
& \begin{array}{cc}
1 \\ r^{(1)} \\ 1
\end{array}\\
(3,0) \longrightarrow \left\{ \begin{array}{ll}
(1,0) \\ (2,0) \\ (3,1) \\ (3,2)
\end{array}\right.
& \begin{array}{cc}
1 \\ 1 \\ \lambda_1^{(2)}n_1(x,\eta^{(2)}) + \lambda_2^{(2)} n_3(x,\eta^{(2)}) \\ r^{(2)}
\end{array}\\
(3,2) \longrightarrow \left\{ \begin{array}{ll}
(3,3) \\ (1,0) \\ (2,2)
\end{array}\right.
& \begin{array}{cc}
\lambda_1^{(2)}n_1(x,\eta^{(2)}) + \lambda_2^{(2)} n_3(x,\eta^{(2)}) \\1 \\ 1 
\end{array}\\
\hline
\end{array}
\end{footnotesize}
$
\captionof{table}{Coupled transitions of two symmetric CPRS, via a basic coupling, that do not preserve the partial order. Nevertheless, starting from an initial configuration that does, the coupled dynamics does not reach these states.}
\label{genprop:monotonicity_cns_tab3}
\end{center}

To verify all the rates above are well defined, one decomposes $n_1(x,\eta^{(i)})$ and \\ $n_3(x,\eta^{(i)})$, $(i=1,2)$, as follows
\begin{align*}
& \begin{multlined} n_1(x,\eta^{(2)}) = |\{ y \sim x : \eta^{(2)}(y)=\eta^{(1)}(y)=1\}| \\
 + |\{ y \sim x : \eta^{(2)}(y)=1, \eta^{(1)}(y)=3\}| + |\{ y \sim x : \eta^{(2)}(y)=1, \eta^{(1)}(y) \in \{0,2\} \}|, \end{multlined}\\
& n_3(x,\eta^{(2)}) = |\{ y \sim x : \eta^{(2)}(y)=\eta^{(1)}(y)=3\}| 
+ |\{ y \sim x : \eta^{(2)}(y)=3, \eta^{(1)}(y)\in \{0,2\}\}|,\\
& n_1(x,\eta^{(1)}) = |\{ y \sim x : \eta^{(2)}(y)=\eta^{(1)}(y)=1\}| \\
& \begin{multlined}n_3(x,\eta^{(1)}) = |\{ y \sim x : \eta^{(2)}(y)=\eta^{(1)}(y)=3\}| + |\{ y \sim x : \eta^{(2)}(y)=1, \eta^{(1)}(y)=3\}| \end{multlined},
\end{align*}
in which case, we decompose the rate 
\begin{small}
\begin{align}\label{genprop:monotonicity_cns_3}
& \lambda_1^{(2)}n_1(x,\eta^{(2)})- \lambda_1^{(1)} n_1(x,\eta^{(1)}) + \lambda_2^{(2)} n_3(x,\eta^{(2)}) - \lambda_2^{(1)}n_3(x,\eta^{(1)}) \nonumber\\
& \quad = (\lambda_1^{(2)}- \lambda_1^{(1)}) |\{ y \sim x : \eta^{(2)}(y)=\eta^{(1)}(y)=1\}| + (\lambda_1^{(2)}- \lambda_2^{(1)}) |\{ y \sim x : \eta^{(2)}(y)=1,\eta^{(1)}(y)=3\}| \nonumber\\ 
& \qquad + (\lambda_2^{(2)}- \lambda_2^{(1)}) |\{ y \sim x : \eta^{(2)}(y)=\eta^{(1)}(y)=3\}| + \lambda_1^{(2)} |\{ y \sim x : \eta^{(2)}(y)=1,\eta^{(1)}(y)\in \{0,2\}\}| \nonumber\\
& \qquad + \lambda_2^{(2)} |\{ y \sim x : \eta^{(2)}(y)=3,\eta^{(1)}(y)\in \{0,2\}\}|
\end{align}
\end{small}
which is non-negative under conditions 1. to 4. coming from (I) and (III) in inequalities  \eqref{genprop:monotonicity_cns1}-\eqref{genprop:monotonicity_cns2}.

\bigskip
Rates of Table \ref{genprop:monotonicity_cns_tab2} are non-negative thanks to conditions 6. to 8., given by inequalities (II)-(i)(ii) with $ \beta=B, \delta=C$. Condition 5. is used by Tables \ref{genprop:monotonicity_cns_tab1} and \ref{genprop:monotonicity_cns_tab3} that correspond to a basic coupling while Table \ref{genprop:monotonicity_cns_tab2} uses a different coupling.  Table \ref{genprop:monotonicity_cns_tab3} listing transitions of the coupled process starting from configurations that do not preserve the defined partial order, nevertheless, starting from an initial configuration that does, dynamics of the coupling given by Tables \ref{genprop:monotonicity_cns_tab1} and \ref{genprop:monotonicity_cns_tab2} do not reach states of Table \ref{genprop:monotonicity_cns_tab3} so that marginals are retrieved.

For a coupled process $(\eta_t^{(1)},\eta_t^{(2)})_{t \geq 0}$ starting from an initial configuration such that $\eta_0^{(1)} \leq \eta_0^{(2)} $, since transitions of the two first Tables preserve the order on $F$, the markovian coupling we just constructed is increasing:
\begin{equation}
\widetilde{\mathbb{P}}^{(\eta^{(1)}_0,\eta^{(2)}_0)} ( \eta^{(1)}_t \leq \eta^{(2)}_t ) = 1 \mbox{ for all } t > 0
\end{equation}
where $ \widetilde{\mathbb P}^{(\eta^{(1)}_0,\eta^{(2)}_0)} $ stands for the distribution of $ (\eta^{(1)}_t, \eta^{(2)}_t)_{t \geq 0} $ starting from $ (\eta^{(1)}_0, \eta^{(2)}_0)$.

We now investigate sufficient conditions for attractiveness, 
under which we build a coupled process through a basic coupling as follows.

\begin{proposition}\label{genprop:monotonicity_cs}
The symmetric process $(\eta_t)_{t \geq 0}$ is monotone, in the sense that one can construct on a same probability space two symmetric processes $(\eta_t^{(1)})_{t \geq 0}$ and $(\eta_t^{(2)})_{t \geq 0}$ with respective parameters $(\lambda_1^{(1)}, \lambda_2^{(1)},r^{(1)})$ and$(\lambda_1^{(2)}, \lambda_2^{(2)},r^{(2)})$ satisfying $\eta_0^{(1)}, \eta_0^{(2)} \in \{0,1\}^{\mathbb Z^d}$, such that 
\begin{equation}
\eta^{(1)}_0 \leq \eta^{(2)}_0 \Longrightarrow \eta^{(1)}_t \leq \eta^{(2)}_t \mbox{ for all } t \geq 0 \mbox{ a.s.}
\end{equation}
if the parameters satisfy
\begin{multicols}{3}
\begin{enumerate}[label={\arabic*.}]
  \item $\lambda_2^{(1)} \leq \lambda_1^{(1)}$, 
    \item $\lambda_2^{(2)} \leq \lambda_1^{(2)}$,
\item $\lambda_1^{(1)} \leq \lambda_1^{(2)}$, 
\item $\lambda_2^{(1)} \leq \lambda_2^{(2)}$,
\item $r^{(1)} \geq r^{(2)}$
\end{enumerate}
\end{multicols}
\end{proposition}
\begin{rmk}\label{genprop:monotonicity_rmk}
In view of the proof of Proposition \ref{genprop:monotonicity_cs}, one can actually relax the admissible initial conditions: it is enough to assume that $\eta_0^{(1)}$ and $\eta_0^{(2)}$ satisfy $\eta_0^{(1)} \leq \eta_0^{(2)}$ and for all $x \in \mathbb Z^d$, $(\eta_0^{(1)}(x),\eta_0^{(2)}(x)) \neq (0,3)$. In particular one could start from $\eta_0^{(1)} = \eta_0^{(2)}$.
\end{rmk}

\begin{proof}
Given initial conditions, possible states for the coupled process keep laying in Table \ref{genprop:monotonicity_cns_tab1} of Proposition \ref{genprop:monotonicity_cns} and the coupled process  does not reach any state of Tables \ref{genprop:monotonicity_cns_tab2} and \ref{genprop:monotonicity_cns_tab3}. One can therefore omit conditions 4. to 6. of the previous Proposition \ref{genprop:monotonicity_cns} and transition rates from the couple $(0,3)$ can be defined through a basic coupling even if it does not preserve the order:
\begin{center}
$
\begin{footnotesize}
\begin{array}{|c|c|}
\hline
\mbox{transition} & \mbox {rate}\\
   \hline
(0,3) \longrightarrow \left\{ \begin{array}{ll}
(1,3)\\(0,1)\\(0,2)\\(2,3)
\end{array}\right.
& \begin{array}{cc}
\lambda_1^{(1)}n_1(x,\eta^{(1)}) + \lambda_2^{(1)} n_3(x,\eta^{(1)}) \\  1 \\ 1 \\ r^{(1)}
\end{array}\\
\hline
\end{array}
\end{footnotesize}
$
\captionof{table}{}\label{genprop:monotonicity_cs_tab2}
\end{center}
in which case, Table \ref{genprop:monotonicity_cs_tab2} substitutes Table \ref{genprop:monotonicity_cns_tab2}.
\end{proof}

If $(\eta_t^{(1)})_{t \geq 0}$ and $(\eta_t^{(2)})_{t \geq 0}$ differ by at most one parameter, one deduces \begin{corollary}\label{genprop:monotonicity_cs_cor} Suppose $\eta_0^{(1)}, \eta_0^{(2)} \in \{0,1\}^{\mathbb Z^d}$. Then for the processes $(\eta_t^{(1)})_{t \geq 0}$ and $(\eta_t^{(2)})_{t \geq 0}$ with parameters $(\lambda_1^{(1)}, \lambda_2^{(1)}, r^{(1)})$ and $(\lambda_1^{(2)}, \lambda_2^{(2)}, r^{(2)})$ respectively, one has
\begin{enumerate}[label=(\roman{*})]
\item Attractiveness: if $(\lambda_1^{(1)},\lambda_2^{(1)},r^{(1)}) = (\lambda_1^{(2)}, \lambda_2^{(2)},r^{(2)})$ , then $\eta_0^{(1)} \leq \eta_0^{(2)} \Rightarrow \eta_t^{(1)} \leq \eta_t^{(2)}$ a.s., for all $t \geq 0$.
\item Increase w.r.t. $\lambda_1$: if $(\eta_t^{(1)})_{t \geq 0}$ and $(\eta_t^{(2)})_{t \geq 0}$ have respective parameters $(\lambda_1^{(1)}, \lambda_2, r)$ and $(\lambda_1^{(2)}, \lambda_2, r)$ such that $\lambda_2 \leq \lambda_1^{(1)} \leq \lambda_1^{(2)}$, then $\eta_0^{(1)} \leq \eta_0^{(2)} \Rightarrow \eta_t^{(1)} \leq \eta_t^{(2)}$ a.s., for all $t \geq 0$.
\item Increase w.r.t. $\lambda_2$: if $(\eta_t^{(1)})_{t \geq 0}$ and $(\eta_t^{(2)})_{t \geq 0}$ have respective parameters $(\lambda_1, \lambda_2^{(1)}, r)$ and $(\lambda_1, \lambda_2^{(2)}, r)$ such that $\lambda_2^{(1)} \leq \lambda_2^{(2)} \leq \lambda_1$, then $\eta_0^{(1)} \leq \eta_0^{(2)} \Rightarrow \eta_t^{(1)} \leq \eta_t^{(2)}$ a.s., for all $t \geq 0$.
\item Decrease w.r.t. $r$: if $(\eta_t^{(1)})_{t \geq 0}$ and $(\eta_t^{(2)})_{t \geq 0}$ have respective parameters $(\lambda_1, \lambda_2, r^{(1)})$ and $(\lambda_1, \lambda_2, r^{(2)})$ such that $r^{(1)} \geq r^{(2)}$ with $\lambda_2 < \lambda_1$, then $\eta_0^{(1)} \leq \eta_0^{(2)} \Rightarrow \eta_t^{(1)} \leq \eta_t^{(2)}$ a.s., for all $t \geq 0$.
\end{enumerate}
\end{corollary}

A highlighted consequence related to Corollary  \eqref{genprop:monotonicity_cs_cor}-(iv) is the monotonicity of the survival probability with respect to the release rate $r$ for fixed $\lambda_1, \lambda_2$:
\begin{corollary}\label{genprop:survivaldecrease}
Suppose $\lambda_2$ and $\lambda_1$ fixed. If $(\eta_t)_{t \geq 0}$ has initial configuration $\eta_0 = \mathbf 1_{\{0\}}$, the mapping  $$r \longmapsto \mathbb P_r (\forall t \geq 0, \ H_t \neq \emptyset) $$ is a non-increasing function.
\end{corollary}
\begin{proof}
Let $(\eta_t^{(1)})_{t \geq 0}$ and $(\eta_t^{(2)})_{t \geq 0}$ be two processes such that $\eta_0^{(1)}, \eta_0^{(2)} \in \{0,1\}^{\mathbb Z^d}$, with respective parameters $(\lambda_1, \lambda_2, r^{(1)})$ and $(\lambda_1, \lambda_2, r^{(2)})$ such that $r^{(1)} \leq r^{(2)}$, then by Corollary \ref{genprop:monotonicity_cs_cor}, for all $t \geq 0$, $$H_0^{(2)} \subset H_0^{(1)} \Longrightarrow H_t^{(2)} \subset H_t^{(1)}.$$
\end{proof}
Analogously to Propositions \ref{genprop:monotonicity_cns} and \ref{genprop:monotonicity_cs}, we can obtain conditions for the asymmetric CPRS to be attractive. For the rest of the section, we refer the reader to \cite{k} for more detailed proofs.
\begin{proposition}\label{genprop:monotonicity_cs_as}
The asymmetric CPRS $(\eta_t)_{t \geq 0}$ is monotone, in the sense that conditions 
\begin{multicols}{3}
\begin{enumerate}[label={\arabic*.}]
  \item $\lambda_2^{(1)} \leq \lambda_1^{(1)}$, 
    \item $\lambda_2^{(2)} \leq \lambda_1^{(2)}$,
\item  $\lambda_1^{(1)} \leq \lambda_1^{(2)}$,
\item $\lambda_2^{(1)} \leq \lambda_2^{(2)}$,
\item $r^{(1)} \geq r^{(2)}$, 
\end{enumerate}
\end{multicols}
\noindent are sufficient to construct on a same probability space two asymmetric CPRS $(\eta_t^{(1)})_{t \geq 0}$ and $(\eta_t^{(2)})_{t \geq 0}$ with respective parameters $(\lambda_1^{(1)}, \lambda_2^{(1)},r^{(1)})$ and $(\lambda_1^{(2)}, \lambda_2^{(2)},r^{(2)})$ and with initial condition $\eta_0^{(1)}, \eta_0^{(2)} \in \{0,1\}^{\mathbb Z^d}$, such that
\begin{equation}
\eta^{(1)}_0 \leq \eta^{(2)}_0 \Longrightarrow \eta^{(1)}_t \leq \eta^{(2)}_t \mbox{ a.s.,  for all } t \geq 0.
\end{equation}
\end{proposition}
\begin{proof}
As in the proof of Proposition \ref{genprop:monotonicity_cs}, one applies \cite[Theorem 2.4]{Bor11} with $j_1,h_1 \in \{0,1\}$ to two asymmetric processes  $(\eta_t^{(1)})_{t \geq 0}$ and $(\eta_t^{(2)})_{t \geq 0}$ with respective parameters $(\lambda_1^{(1)}, \lambda_2^{(1)},r^{(1)})$ and $(\lambda_1^{(2)}, \lambda_2^{(2)},r^{(2)})$. Using 
\eqref{genprop:monotonicity_cns1}-\eqref{genprop:monotonicity_cns2} with rates given by \eqref{borrello_CPDRE}, 
we deduce the following necessary conditions from the first inequality:
\begin{enumerate}[label={(\Roman*)}]
\item $j_1 \in \{0,1\}$, $\delta=\beta=B$ in \eqref{genprop:monotonicity_cns1} give
\begin{enumerate}[label={(\roman*)}]
\item $\alpha=\gamma=C$, $\beta=B$, $\delta=C$: $\lambda_2^{(1)} \leq \lambda_2^{(2)}$ stated by condition 4.
\item $\alpha=C, \gamma=D$: $\lambda_2^{(1)} \leq \lambda_1^{(2)}$ stated by conditions 1. and 3.
\item $\alpha=\gamma=D$: $\lambda_1^{(1)} \leq \lambda_1^{(2)}$ stated by condition 3.
\end{enumerate}
\item $j_1=0$, $\beta=B$, $\delta=1+\beta=C$ in \eqref{genprop:monotonicity_cns1} give
\begin{enumerate}[label={(\roman*)}]
\item $\alpha=D$: $\lambda_1^{(1)} \leq 1$.
\item$\alpha=C$: $\lambda_2^{(1)} \leq 1$.
\end{enumerate}
\end{enumerate}
The relation \eqref{genprop:monotonicity_cns2}  staying unchanged, one has
\begin{enumerate}[label=(\Roman*)]
  \setcounter{enumi}{2}
\item $h_1=0$, $\gamma=\alpha \in \{B,D\}$ in \eqref{genprop:monotonicity_cns2} give $r^{(1)} \geq r^{(2)}$ stated by condition 5.
\item $h_1=0$,$\alpha=B,\gamma=1+\alpha=C$ in \eqref{genprop:monotonicity_cns2} give $r^{(1)} \geq  1$.
\end{enumerate}
To sum up, the obtained necessary conditions are 
\begin{multicols}{4}
\begin{enumerate}[label={\arabic*.}]
  \item $\lambda_2^{(1)} \leq \lambda_1^{(1)}$, 
    \item $\lambda_2^{(2)} \leq \lambda_1^{(2)}$,
\item $\lambda_1^{(1)} \leq \lambda_1^{(2)}$, 
\item $\lambda_2^{(1)} \leq \lambda_2^{(2)}$,
\item $r^{(1)} \geq r^{(2)}$
\item $\lambda_1^{(1)} \leq 1$, 
\item $\lambda_2^{(1)}  \leq 1$, 
\item $r^{(1)} \geq 1.$
\end{enumerate}
\end{multicols}
As for Proposition \ref{genprop:monotonicity_cns}, these conditions allow us to construct an increasing markovian coupling and are thus necessary and sufficient. But as in Proposition \ref{genprop:monotonicity_cs}, given our initial configurations, state $(0,3)$ is not possible for the coupled process. One can thus dispense conditions 6 to 8. and conditions 1. to 5. are sufficient to construct a basic coupling.
\end{proof}

We now compare asymmetric, symmetric and basic contact processes.
\begin{proposition}\label{genprop:monotonicity_sym_higher}
Let $(\eta_t)_{t \geq 0}$ be an asymmetric CPRS and $(\chi_t)_{t \geq 0}$ be a symmetric CPRS, both with parameters  $(\lambda_1, \lambda_2,r)$ and $\eta_0, \chi_0 \in \{0,1\}^{\mathbb Z^d}$ such that $\lambda_2 < \lambda_1$, then for all $t \geq 0$, $$\eta_0 \leq \chi_0 \Rightarrow \eta_t \leq \chi_t \mbox{ a.s.  for all } t \geq 0$$
\end{proposition}
\begin{proof}
Consider a symmetric process with rates \eqref{borrello_CPDRE}-\eqref{borrello_CPDRE2} and an asymmetric process with rates \eqref{borrello_CPDRE}. Similarly to the proof of Proposition \ref{genprop:monotonicity_cns}, we apply inequalities \eqref{genprop:monotonicity_cns1}-\eqref{genprop:monotonicity_cns2} to these rates so that they give the following necessary conditions:
\begin{enumerate}[label={(\Roman*)}]
\item $j_1 \in \{0,1\}$, $\delta=\beta=B$,$\alpha=C,\gamma=D$ in \eqref{genprop:monotonicity_cns1} give: $\lambda_2 \leq \lambda_1$
\item $h_1=0$, $\alpha=B,\gamma=1+\alpha=C$  in \eqref{genprop:monotonicity_cns2} give $r \geq  1$
\end{enumerate}
As previously, condition $r \geq 1$ is necessary to construct an increasing markovian coupled process in a general framework, but if one restricts the initial conditions to satisfy $\eta_0 \leq \chi_0$ and $\eta_0,\chi_0 \in \{0,1\}^{\mathbb Z^d}$, condition $\lambda_2 < \lambda_1$ is sufficient and the coupled process can be constructed through a basic coupling.
\end{proof}


\begin{proposition}\label{genprop:monotonicity_cp_higher}
Let $(\xi_t)_{t \geq 0}$ be a basic contact process on $\{0,1\}^{\mathbb Z^d}$ with growth rate $\lambda_1$ and let $(\chi_t)_{t \geq 0}$ be a symmetric CPRS with parameters  $(\lambda_1,\lambda_2,r)$ such that $\lambda_2 < \lambda_1$. Then,
$$\chi_0 \leq \xi_0 \Rightarrow \chi_t \leq \xi_t \mbox{ a.s. for all } t \geq 0$$
\end{proposition}
\begin{proof}
Consider a symmetric CPRS with rates \eqref{borrello_CPDRE}-\eqref{borrello_CPDRE2} and a supercritical contact process $(\xi_t)_{t \geq 0}$ viewed as a symmetric CPRS with parameters  $(\lambda_1^{(2)},\lambda_2^{(2)},r^{(2)})$ with $\lambda_1^{(2)} = \lambda_1$, $\lambda_2^{(2)} = 0$, $r^{(2)}=0$, that is with rates \eqref{borrello_contactsuper}.
Applying \eqref{genprop:monotonicity_cns1}-\eqref{genprop:monotonicity_cns2} to rates \eqref{borrello_CPDRE}-\eqref{borrello_CPDRE2} and \eqref{borrello_contactsuper}, 
since values $A$ and $C$ do not exist for the process $\xi_t$, 
we have the following necessary condition:
$j_1 \in \{0,1\}$, $\beta=\delta=B$, $\alpha=C$, $\gamma=D$ in \eqref{genprop:monotonicity_cns1} give $\lambda_2 \leq \lambda_1$. While relation  \eqref{genprop:monotonicity_cns2} does not give further condition. Condition $\lambda_2 \leq \lambda_1$ is sufficient and allows us to construct a coupled process through basic coupling.
\end{proof}
%
%
\begin{proposition}\label{genprop:monotonicity_cp_lower}
Let $(\eta_t)_{t \geq 0}$ be a symmetric CPRS with parameters $(\lambda_1,\lambda_2,r)$ and $(\widetilde \xi_t)_{t \geq 0}$ a basic contact process on $\{2,3\}^{\mathbb Z^d}$ with rate $\lambda_2$ such that $\lambda_2 \leq \lambda_1$. Then,
$$\widetilde \xi_0 \leq \eta_0 \Rightarrow \widetilde \xi_t \leq \eta_t \mbox{ a.s., for all } t \geq 0.$$
\end{proposition}
\begin{proof}
For the process $(\widetilde \xi_t)_{t \geq 0}$, values $B$ and $D$ are not reached 
so that the necessary and sufficient conditions on the parameters given by \eqref{genprop:monotonicity_cns1}-\eqref{genprop:monotonicity_cns2}
exhibit the condition:
$j_1 \in \{0,1\}$, $\beta=\delta=A$, $\alpha=C$, $\gamma=D$  in \eqref{genprop:monotonicity_cns1} give $\lambda_2 \leq \lambda_1$. Inequality \eqref{genprop:monotonicity_cns2} gives no condition on the rates and condition $\lambda_2 \leq \lambda_1$ is sufficient to construct the coupled process $(\widetilde \xi_t,\eta_t)_{t \geq 0}$ via a basic coupling.
\end{proof}


\section{Phase transition on $\mathbb Z^d$}\label{sec:phase}
Here, we take advantage of all the stochastic order relations established in Section \ref{sec:order} to derive a phase transition for the CPRS, in both symmetric and asymmetric cases. According to \eqref{cpsp:surviv.extinct}, assume the process to have initial configuration $\eta_0 = \mathbf 1_{\{0\}}$ and note $\eta_t = \eta_t^{\{0\}}$. We first explore cases when $\lambda_2 < \lambda_1$ are both smaller or larger than $\lambda_c$.

\begin{proof}[Proof of Proposition \ref{phase:trivial_sub}]
Let $(\xi_t)_{t \geq 0}$ be a basic contact process with growth rate  $\lambda_1$ and let $(\eta_t)_{t \geq 0}$ be a symmetric CPRS with parameters $(\lambda_1,\lambda_2,r)$ such that $\eta_0 \leq \xi_0$. By Proposition \ref{genprop:monotonicity_cp_higher}, $(\xi_t)_{t \geq 0}$ is stochastically larger than $(\eta_t)_{t \geq 0}$. Since $ \lambda_1 \leq \lambda_c $, $ (\xi_t)_{t \geq 0} $ is subcritical, thus, the symmetric CPRS dies out.

The extinction of the asymmetric CPRS is a consequence of the extinction of the symmetric process and  Proposition \ref{genprop:monotonicity_sym_higher}.
\end{proof}

\begin{proof}[Proof of Proposition \ref{phase:trivial_sup}]
Let $(\widetilde \xi_t)_{t \geq 0}$ be a basic contact process with growth rate $\lambda_2$ on  $\{2,3\}^{\mathbb Z^d}$ and let $(\eta_t)_{t \geq 0}$ be a symmetric CPRS with parameters $(\lambda_1,\lambda_2,r)$. By Proposition \ref{genprop:monotonicity_cp_lower}, $(\xi_t)_{t \geq 0}$ is stochastically lower than $(\eta_t)_{t \geq 0}$. Since $ \lambda_2 > \lambda_c $, the process $ (\widetilde \xi_t)_{t \geq 0} $ is supercritical and therefore, the symmetric CPRS  survives.
\end{proof}

 Before proving Theorem \ref{transition:existence_sym} to \ref{transition:diesout} , the following Subsection \ref{transition:behaviour} deals with consequences of Theorems \ref{transition:existence_sym} and \ref{transition:existence_asym} along with monotonicity results of Section \ref{sec:order}.

\subsection{Behaviour of the critical value with varying growth rates}\label{transition:behaviour}
Here we prove monotonicity between growth rates and the release rate in the sense that, the more competitive the species is (i.e. the higher the parameter $\lambda_2$ is) or the fitter the species is (i.e. the higher the parameter $\lambda_1$ is), the higher the release rate is (i.e. the higher the critical value $r_c$ is).
\begin{proposition} For $j=1,2,$ the function $\lambda_j \longmapsto r_c (\lambda_j)$ is non-increasing.
\end{proposition}
\begin{proof}
We consider $j=2$, the case $j=1$ is similar. Let $(\eta_t)_{t \geq 0}$ and $(\eta'_t)_{t \geq 0}$ be two CPRS with respective parameters $(\lambda_1, \lambda_2, r)$ and $(\lambda_1, \lambda'_2, r)$ such that $\lambda_2 < \lambda_2 '$. Existence of $r_c$ and $r_c'$ for those processes is given by Theorems \ref{transition:existence_sym} and \ref{transition:existence_asym}. 
By contradiction, suppose $r_c > r_c'$. Let $r$ be fixed be such that $r_c > r > r_c'$. 
By Corollary \ref{genprop:monotonicity_cs_cor}-(iii), if $\eta_0 = \eta'_0$ then $\eta_t \leq \eta'_t$ a.s. and by Theorem \ref{transition:diesout} and Corollary \ref{genprop:monotonicity_cs_cor} (iv),
$$\mathbb P_r(\forall t \geq 0, \ H'_t \neq \emptyset) \leq \mathbb P_{r'_c} (\forall t \geq 0, \ H_t' \neq \emptyset) = 0$$
But since $r<r_c$, the process $(\eta_t)_{t \geq 0}$ survives: $\mathbb P_r(\forall t \geq 0, \ H_t \neq \emptyset) > 0$. This contradicts $\eta_t \leq \eta'_t$ a.s., hence $r_c \leq r'_c $.
\end{proof}

\subsection{Subcritical case}\label{transition:subcritical}
The following proof relies on a comparison with an oriented percolation process on the grid $\mathcal{L}= \{ (x,n) \in \mathbb Z^2 : x+n \mbox{ is even}, n \geq 0 \}$. We follow arguments used in \cite{KST04} for a spatial epidemic model.

\begin{proof}[Proof of Theorem \ref{transition:existence_sym} (i)]
To simplify notations, choose $d=1$ but the proof remains the same for any $d \geq 2$. Introduce the space-time regions, 
\begin{align*}
& \mathcal{B} = (-4L, 4L) \times [0,T], \ \ \mathcal{B}_{m,n} = (2m L e_1 , nT) + \mathcal{B}  \\
& I = [-L, L], \ \ I_m = 2mLe_1 + I
\end{align*}
for positive integers $L, T$ to be chosen later, where $(e_1,...,e_d)$ denotes the canonical basis of $\mathbb R^d$. 

\bigskip
Consider the process $(\eta_t^{m,n})_{t \geq 0}$ restricted to the region $\mathcal{B}_{m,n}$, that is, constructed from the graphical representation where only arrival times of the Poisson processes occurring within $\mathcal{B}_{m,n}$ are taken into account. Therefore, a birth on a site $x \in \mathcal B_{m,n}$ from some site $y$ only occurs if $y \in \mathcal B_{m,n}$. By Proposition \ref{genprop:monotonicity_cs} and Remark \ref{genprop:monotonicity_rmk}, one has 
\begin{equation}\label{subcritical:order}
\eta_t^{m,n} \leq \eta_t \big|_{\mathcal B_{m,n}},
\end{equation}
for all $t > 0$ if $\eta_0^{m,n} =\eta_0 \big|_{\mathcal B_{m,n}}$.

Let $k= \lfloor \sqrt L \rfloor$ and define $\mathcal C = [-k,k]$. One declares $ (m,n) \in \mathcal{L} $ to be \textsl{wet} if for any configuration at time $nT$ such that there is a translate of $\mathcal C$ full with 1's in $I_m$ with $I_m$ containing only 0's and 1's, the process restricted to $\mathcal B_{m,n}$ is such that at time $(n+1)T$ there are a translate of $\mathcal C$ in $I_{m-1}$ and a translate of $\mathcal C$ in $I_{m+1}$, both full of 1's, with $I_{m-1}$ and $I_{m+1}$ containing only 0's and 1's.

\bigskip
We show the probability of a site $(m,n) \in \mathcal L$ to be wet can be made arbitrarily close to 1 for $L$ and $T$ chosen sufficiently large. By translation invariance, it is enough to deal with the case $(m,n) = (0,0)$.

\bigskip
Suppose $I$ contains only 0's and 1's as well as the translate of $\mathcal C$ full of 1's and set $r=0$ in $\mathcal B$, that is, no more type-2 individuals arrive in the box  $\mathcal B$ after time 0. If type-2 individuals are present on $(-4L,-L)\cup (L,4L) \times \{0\}$, the probability of the event $E$ ``they all die by time $T/2$'' is at least 
$$\Big( 1- \exp (-T/2) \Big)^{6L}$$ 
which is larger than $1-\epsilon$ for $T$ and $L$ large enough. On $E$, the process restricted to the box $\mathcal B$ is now from time $T/2$ a supercritical contact process $(\xi^{m,n}_t)_{t \geq T/2}$ with distribution $\widetilde{\mathbb P}(\xi^{m,n}_t \in \cdot)$. 
Define $\tau (\ell)$, the hitting time of the trap state $\underline 0$ of the contact process starting from $[-\ell,\ell]$ and restricted to $[-\ell,\ell] \times [0,T/2]$. To ensure there are still enough 1's for $\xi_{T/2}^{m,n}$, T. Mountford \cite{Tom93} has proven
\begin{equation}\label{mountford93}
\widetilde{\mathbb P} \big(\tau(\ell) \leq \exp(\ell) \big) \leq \exp(-\ell) \ \mbox{ for } \ell \mbox{ large enough}
\end{equation}
Partition $\mathcal C$ into $M= \lfloor \sqrt{k} \rfloor$ boxes, each of them being a translate of $[0,M]$. From each of these $M$ boxes, say box $j \leq M$, run a supercritical contact process $(\zeta_t^j)_{t \geq 0}$ which coincides with the restriction of $\xi_t^{m,n}$ to this box. 
Then for the union of these $j$ boxes 
the probability there is at least  $M$ 1's within $\mathcal C$ by time $T/2$ is after \eqref{mountford93}, with $T$ such that $\exp(M) \leq T/2$, at least
\begin{equation}
\widetilde{\mathbb P} (\tau(M) \geq T/2 )^M \geq \widetilde{\mathbb P} (\tau(M) \geq \exp(M) )^M \geq (1-\exp(-M) )^M
\end{equation}
which can be made larger than $1-\epsilon$, for $M$, i.e. $L$, large enough. 

\bigskip
R. Durrett and R. Schinazi \cite{DS93} have shown that for a contact process $(\xi_t)_{t \geq 0}$, for any $A \subset \mathbb Z$, except for a set with exponentially small probability, either $\Xi_t^A = \emptyset$, or $\xi_t^A=\xi_t^{\mathbb Z}$ on a linearly time growing set $[-\alpha t, \alpha t]$:  
there exists $\alpha > 0$ such that for all $A \subset \mathbb Z$, there exist constants $C, \gamma$ such that for  $x \in A + \alpha t$,
\begin{equation}\label{ds}
\widetilde{\mathbb P} \big(\Xi_t^A \neq \emptyset, \ \xi_t^A(x) \neq \xi_t^{\mathbb Z} (x)\big) \leq C\exp (-\gamma t )
\end{equation}

We applied this result with $A \subset \mathcal C$, the set of 1's in the box. We just proved that $|A| > |M|$. Moreover by attractiveness (i.e. Proposition \ref{genprop:monotonicity_cns}), one can choose $k$, and so $L$, large enough so that this supercritical contact process $\xi_t^A$ starting from at least $M$ 1's survives at time $T/2$ with probability close to 1, hence $\xi_t^A \neq \emptyset$ and \eqref{ds} is valid. In this situation, taking $T/2 = 9L/(2 \alpha )$ with $L$ large enough, the process $(\xi_t^{m,n})_{t \geq 0}$ starting from at least  $M$ 1's in $[-L,L]$ at time $T/2$ will be coupled with a process $\xi_t^{\mathbb Z}$ on $[-3L,3L]$ with probability at least $1-\epsilon$ at time $T$. Hence, since $3L > \alpha T /2 > 2L$, by time $T$ the contact process $\xi_t^A$ started inside $[-L,L]$ has not reached the boundary of $[-4L,4L]$ with probability close to 1. Then, the process $\xi_t^{m,n}$ and the contact process $\xi_t^A$ are the same with probability $1-\epsilon$ in $[-4L,4L]$ ; this way, the coupling of $(\xi_t^{m,n})_{t \geq 0}$ with $(\xi_t^{\mathbb Z})_{t\geq 0}$ works so far with probability $1-\epsilon$ if $L$ is large enough.

\bigskip
Since the distribution of $(\xi_t^{\mathbb Z})_{t \geq 0}$ is stochastically larger than the upper invariant measure $\overline{\nu}$ of the contact process, on the survival event, $\overline{\nu}$ loads a positive density $\rho$ of 1's. Since $\overline{\nu}$ is ergodic,
$$\lim\limits_{ L \rightarrow \infty} \dfrac{1}{2L+1} \sum\limits_{x =-3L}^{-L} \mathbf{1}\{\eta(x)=1\} = \rho \ \ \overline{\nu} \mbox{ -a.e.}$$
Then if $L$ is large enough, there are at least $k$ 1's in any interval of length $2L+1$ with $ \overline{\nu}$-probability at least $1-\epsilon$. 
Since we obtained that $(\xi_t^{m,n})_{t \geq 0}$ is coupled to $(\xi_t^{\mathbb Z})_{t \geq 0}$ by time $T$ with probability at least $1-2\epsilon$, for $L$ large enough, there are at least $k$ 1's in $[-3L,-L]$ at time $T$  with probability at least $1-2\epsilon$ and similarly, at least $k$ 1's at time $t$ in $[L,3L]$ with probability at least $1-2\epsilon$ as well for $(\xi_t^{m,n})_{t \geq 0}$. Consequently,
\begin{align}\label{wetsurvie}
\widetilde{\mathbb P} ((0,0) \mbox{ wet} ) > 1 - 4\epsilon, \mbox{ if } r = 0.
\end{align}

Since $\mathcal{B}$ is a finite space-time region, for fixed $L,T$, one can pick $r_0>0$ small enough so that if $r \leq r_0$
$$\mathbb P_r ((0,0) \mbox{ wet} ) > 1 - 6\epsilon.$$

Now construct a percolation process by defining $G_{m,n}=\{ (m,n) \mbox{ wet} \}$. Note $G_{m,n}$ depends only on the process constructed in $\mathcal{B}_{m,n}$, and for $(a,b) \in \mathcal{L}$, events $G_{m,n}$ and $G_{a,b}$ are independent if $(m,n)$ and $(a,b)$ are not neighbours. The events $\{G_{m,n}, (m,n) \in \mathcal{L}\}$ are thus 1-dependent. By the comparison theorem \cite[Theorem 4.3]{Dur93}, the process $(\eta_t^{m,n})_{t \geq 0}$ restricted to regions $\mathcal{B}_{m,n}$ is stochastically larger than a 1-dependent percolation process with probability $1-\epsilon$.

\bigskip
One can choose $\epsilon$ small enough so that percolation occurs in the 1-dependent percolation process with density $1-\epsilon$, see \cite[Theorem 4.3]{Dur93}.
\end{proof}

\subsection{Supercritical case}\label{transition:supercritical}

In the following, one compares our particle system with a percolation process on $\mathbb Z^2 \times \mathbb Z_+$ using arguments exposed as in \cite{SS06}.

\begin{proof}[Proof of Theorem \ref{transition:existence_sym} (ii)]
Assume $d =2$, the proof can similarly be extended to higher dimensions. For all $(k,m,n) \in \mathbb Z^2 \times \mathbb Z_+$.  Introduce the following space-time regions, for positive $L$ and $T$ to be chosen later.
$$\begin{array}{lll}
& \mathcal{A} = [-2L, 2L]^2 \times [0, 2T] & \hspace*{3cm} \mathcal{A}_{k,m,n} = \mathcal{A} + (kL, mL, nT) \\
& \mathcal{B} =  [-L, L]^2 \times [T, 2T] & \hspace*{3cm} \mathcal{B}_{k,m,n} = \mathcal{B} + (kL, mL, nT) \\
& \mathcal{C} = \mathcal{C}_{bottom} \bigcup \mathcal{C}_{side} & \hspace*{3cm}  \mathcal{C}_{k,m,n} = \mathcal{C} + (kL, mL, nT) \\
\end{array}$$
$\begin{array}{lll}
& \ \mbox{ where } \ \mathcal{C}_{bottom} = \{ (m,n,t) \in \mathcal{A} : t = 0 \} \ \ \mbox{ and }\\
& \qquad \mathcal{C}_{side} = \{ (m,n,t) \in \mathcal{A} : |m|= 2L \mbox{ or } |n|= 2L \} \\
\end{array}$

\bigskip
Consider a restriction of the process $(\eta_t)_{t \geq 0}$ to $\mathcal A_{k,m,n}$, that is, the process $(\eta_t^{k,m,n})_{t \geq 0}$ constructed from its graphical representation within $\mathcal{A}_{k,m,n}$.

One declares a site $(k,m,n) \in \mathbb Z^2 \times \mathbb Z_+$ to be \textsl{wet} if the process $(\eta_t^{k,m,n})_{t \geq 0 }$ contains no wild individual in $\mathcal{B}_{k,m,n}$ starting from any configuration in $\mathcal{C}_{k,m,n}$. Therefore it will be the same for $\eta_t \big|_{\mathcal A_{k,m,n}}$. Sites that are not wet are called \textsl{dry}.\\

For any $\epsilon >0$, we show that for some chosen $L$ and $T$ any site $(k,m,n)$ is wet with probability close to 1 when $r$ is large enough. By translation invariance, it is enough to consider $(k,m,n)=(0,0,0)$. Set $r=\infty$ in $\mathcal A$. Then, the process $(\eta_t^{k,m,n})_{t\geq 0}$ contains only sites in state 2 or 3: sites in state 0 or 1 flip instantaneously  into state 2 and 3 respectively.
That is, $(\eta_t^{k,m,n})_{t \geq 0}$ is in fact a contact process $(\widetilde \xi_t^{k,m,n})_{t \geq 0}$ on $\{2,3\}^{[-2L,2L]}$ The contact process $(\widetilde \xi_t)_{t \geq 0}$ on $\{2,3\}^{\mathbb Z^2}$ with growth rate $\lambda_2<\lambda_c$ is subcritical. 

\bigskip
If there is some wild individual in $\mathcal B$, it should have come from a succession of births started somewhere in $\mathcal C$. Starting from a site in $\mathcal C_{side}$, a path to $\mathcal B$ should last at least $L$ sites ; according to C. Bezuidenhout and G. Grimmett \cite{BG91} there exists such a path with probability at most $C\exp(-\gamma L)$, for some positive constants $C,\gamma$. Hence,
$$\mathbb P ( \exists (x,t) \in \mathcal C_{side} \times [0,2T] : (x,t) \rightarrow \mathcal B ) \leq 4 \big(2T \times (4L+1) \big) C\exp(-\gamma L) $$
Similarly, starting from the base $\mathcal C_{bottom}$, there exists  a path lasting at least $T$ sites with probability 
$$\mathbb P ( \exists (x,t) \in \mathcal C_{bottom} : (x,t) \rightarrow \mathcal B ) \leq  (4L+1)^2 C\exp(-\gamma T) $$
Consequently if $r = \infty$, 
$$\mathbb P((0,0,0) \mbox{ wet } ) \geq 1 - 4 \big(2L \times (4L+1) \big) C e^{-\gamma L} - (4L+1)^2 Ce^{-\gamma T} \geq 1 -\epsilon / 2,$$
for $L$ and $T $ large enough.

Since $ \mathcal A$ is a finite space-time region, one can pick $r$ large enough so that with probability at least $1-\epsilon/2$, an exponential clock with parameters $r$ rings before any other so that there are no type-1 individuals in $\mathcal A$ with probability close to 1:
$$\mathbb P_r (k,m,n) \mbox{ wet} ) \geq 1 - \epsilon$$
for $r$ large enough.

\bigskip
To construct a percolation process on $\mathbb Z^2 \times \mathbb Z_+$, one puts an oriented arrow from $(k,m,n)$ to $(x,y,z)$ if $ n \leq z$ and if  $\mathcal A_{k,m,n} \cap  \mathcal A_{x,y,z} \neq \emptyset$. The event $G_{k,m,n} = \{(k, m, n) \mbox{ wet} \}$ depends only on the graphical construction of the process within $\mathcal A_{k,m,n}$, furthermore, for all $(k, m, n) \in \mathbb Z^2 \times \mathbb Z_+$, there is a finite number of sites $(x,y,z) \in \mathbb Z^2 \times \mathbb Z_+$ such that $\mathcal A_{k,m,n} \cap \mathcal A_{x,y,z} \neq \emptyset$. The percolation process is dependent but of finite range. The existence of a path of wild individuals for the particle system corresponds to a path of dry sites for the percolation and we proved that dry paths are unlikely. This implies that for any site, there will be no more wild individual after a finite random time, see \cite[Section 8 p.16]{BGS98}. By translation invariance, this corresponds to the extinction of the process.
%
%
%
%
\end{proof}


\bigskip
This proves the existence of a phase transition for the symmetric CPRS. The proof of Theorem \ref{transition:existence_sym}-(i) only uses that contact process with growth rate $\lambda_1$ is supercritical, hence the existence of $s_0$ in Theorem \ref{transition:existence_asym}-(i) as well. By Proposition \ref{genprop:monotonicity_sym_higher}, the asymmetric CPRS dies out if the symmetric one does so that existence of $s_1$ in Theorem \ref{transition:existence_asym}-(ii) follows from Theorem \ref{transition:existence_sym}-(ii).  Though, one remarks conditions of Theorem \ref{transition:existence_asym} are milder: one can actually show the existence of $s_1$ in a neater way, assume $\lambda_2 > \lambda_c$: without  transition "$2 \rightarrow 3$" in the asymmetric case and choosing $r=\infty$, notice for a subcritical contact process on $\{2,3\}^{\mathbb Z^2}$, there are no paths of wild individuals created by the 3's from the boundary $\mathcal C_{k,m,n}$ up to extinction, but this occurs exponentially fast, see \cite{BG91}.

\subsection{Critical case}


Note the arguments developed by \cite{BG90} rely on elementary properties of the contact process making them robust. They are also valid for the CPRS since it satisfies the following properties.
\begin{enumerate}[label=(\Alph*)]
\item \label{critical:H1} contact process-like dynamics: one retrieves the growth rate $\lambda_1$ or $\lambda_2$ of a basic contact process, even if it is determined randomly. 
\item \label{critical:H2} attractiveness, cf. Section \ref{sec:order}.
\item \label{critical:H3} correlation inequalities: such as FKG inequality.
\end{enumerate}
Note that the use of \ref{critical:H3} is possible because the proof works in finite space-time regions and since survival events are monotone functions of Poisson processes (from the graphical representation, see Section \ref{sec:graph}), correlation inequalities for such Poisson processes apply, see \cite[Lemma 2.1]{Jan84} for instance. 
We therefore omit the proof of Theorem \ref{transition:diesout} but rather present a brief outline and refer the reader to \cite{BG90,Lig99,SW08} for further details.

\bigskip
\textit{Outline of the proof of Theorem \ref{transition:diesout}}. The first step consists to observe that if the process survives in an arbitrary large box, then, that it reaches its boundaries densely. We shall estimate these densities at each side of a space-time region.

This way, one can repeat this step by running the process in an other adjacent box starting from the boundary of the previous one and so on, conditionally on the fact that the starting configuration is dense enough. This is the second step. In connection with the proofs of Theorems \ref{transition:existence_sym} and \ref{transition:existence_asym} where we looked after having translations of occupied finite intervals at a given fixed time, here we look after having translations of the densities in some space-time slab.

Now, knowing that at each stage, one can construct overall a path of adjacent boxes wherein the process survives and reaches the boundaries densely, it remains to compare the process with an oriented percolation process to extend the arguments to infinite space and time. As before, compare a space-time box to a vertex in the even lattice of an oriented percolation so that one declare a vertex to be wet if some good event associated to the box is a success. Conclude thanks to well-settled results about percolation theory.

\bigskip
To conclude the section, Theorem \ref{transition:sym} (resp. Theorem \ref{transition:asym})
follows from 
Theorem \ref{transition:existence_sym} (resp. Theorem \ref{transition:existence_asym}) and by monotonicity arguments given by Corollary \ref{genprop:survivaldecrease}.

\section{Mean-field model}\label{sec:mean}
The mean-field model associated to the CPRS is a deterministic and non-spatial model where all individuals are mixed up, featuring the densities of each type of individuals. Such models give rise to differential systems and are interesting to compare with stochastic behaviours. We investigate here the equilibria, in both asymmetric and symmetric models, since existence of such equilibria yield the existence of a critical value. We exhibit conditions on $r$ to deduce bounds on the critical value.

\bigskip
Subsequently, let $u_i$ be the density of type-$i$ individuals for $i=1,2,3$. Overall, one has $u_1+u_2+u_3 = 1-u_0$. Furthermore, in connection with the definition of wild and sterile individuals, we consider as well $v_1$, resp. $v_2$, the density of the wild individuals (type-1 and type-3 individuals), resp. the sterile individuals (type-2 and type-3 individuals), and the density of empty sites $v_0 = u_0$.
Relations between the $u$-system and the $v$-system are described by
\begin{equation}\label{u-v}
\left\{
\begin{array}{lll}
u_1 = 1-v_0-v_2 \\
u_2 = 1-v_0-v_1 \\
u_3 = v_0+v_1+v_2-1
\end{array}
\right. .
\end{equation}
Since we consider densities, both systems satisfy
\begin{equation}\label{cond-01}
u_i  \in [0,1]  \mbox{ for } i=0,1,2,3,  \ \ \ \ v_i \in [0,1] \mbox{ for } i=0,1,2.
\end{equation}
Following the most suitable case, we investigate either the $u$-system or both systems.
\subsection{Asymmetric CPRS}
Assuming total mixing, the mean-field model associated to the asymmetric CPRS is given by:

\begin{equation}\label{v_asym}
\left\{
\begin{array}{lll}
v_0' = -2d \Big( (\lambda_2-\lambda_1)v_0 + \lambda_2 v_1 +(\lambda_2-\lambda_1)v_1 + \lambda_1-\lambda_2 \Big) v_0 - (r+2) v_0 - v_1 - v_2 +2\\
v_1' = 2d \Big( (\lambda_2-\lambda_1)v_0 +  \lambda_2 v_1 +   (\lambda_2-\lambda_1)v_2 + \lambda_1-\lambda_2 \Big) v_0 - v_1\\
v_2'= r(1-v_2)-v_2
\end{array}
\right.
\end{equation}
This system gives rise to two equilibria:
$$(v_0,v_1,v_2) = \Big(0,0,\frac{r}{r+1}\Big),$$
$$(v_0,v_1,v_2) = \Big( \frac{2+r}{4d\lambda_1r+2d\lambda_2r}, \frac{r+2}{2(r+1)} - \frac{(r+2)^2}{2(4d\lambda_1+2d\lambda_2r)}, \frac{r}{r+1}\Big). $$
Note that the first equilibrium gives $(u_1,u_2,u_3) = \Big(0,\frac{r}{r+1},0\Big)$, which puts a positive density on the sterile individuals and none on the others, which corresponds to the extinction of the process. By checking conditions \eqref{cond-01} on the second, we get that:  the density $v_1$ is non-negative as soon as
$$4d\lambda_1+2d\lambda_2 r > (r+1)(r+2).$$
which gives the following condition
\begin{equation}\label{cond-r_0}
r > \dfrac{2d\lambda_2 - 3 + \sqrt{(2d\lambda_2-3)^2 -8(1-2d\lambda_1)}}{2}
\end{equation}
This indicates a lower bound for the phase transition.

\subsection{Symmetric CPRS}
The mean-field equations associated to the symmetric CPRS are :
\begin{equation}\label{u_sym}
\left\{
\begin{array}{lll}
u_1' = 2d(\lambda_1 u_1 + \lambda_2 u_3) u_0 + u_3 -(r+1)u_1 \\
u_2' = ru_0 + u_3 - u_2 - 2d(\lambda_1 u_1 + \lambda_2 u_3) u_2 \\
u_3' = ru_1 + 2d(\lambda_1 u_1 + \lambda_2 u_3) u_2 -2u_3
\end{array}
\right.
\end{equation}
As previously, this system admits one trivial equilibrium: 
$$(u_1,u_2,u_3) = \Big(0,\frac{r}{r+1},0\Big)$$
retrieving once again a situation related to the extinction of the process, by a positive density of sterile individuals and none of the wild ones. We derive the non-trivial equilibrium thanks to the corresponding $v$-system:
\begin{equation}\label{v_sym}
\left\{
\begin{array}{lll}
v_0' = -2d \Big( (\lambda_2-\lambda_1)v_0 + \lambda_2 v_1 +(\lambda_2-\lambda_1)v_1 + \lambda_1-\lambda_2 \Big) v_0 - (r+2) v_0 - v_1  v_2 +2\\
v_1' = 2d \Big( (\lambda_2-\lambda_1)v_0 +  \lambda_2 v_1 +   (\lambda_2-\lambda_1)v_2 + \lambda_1-\lambda_2 \Big) )(1-v_1) - v_1\\
v_2'= r(1-v_2)-v_2
\end{array}
\right.
\end{equation}

Let us determine the non-trivial equilibrium. Last line of \eqref{v_sym} gives already $v_2 = \dfrac{r}{r+1}$. Using relations of \eqref{u-v} in \eqref{v_sym}, according to $v_1'= 0$, an equilibrium $(v_0,v_1,v_2)$ satisfies in particular 
\begin{equation}
v_1=2d(\lambda_1u_1 + \lambda_2 u_3)(1-v_1)
\end{equation}
checking $v_1$ cannot be equal to 1, one then has 
\begin{equation}\label{cond-v_1}
\dfrac{v_1}{1-v_1}=2d(\lambda_1u_1 + \lambda_2 u_3)
\end{equation}
and
\begin{equation}\label{v_1}
v_1 \neq 1.
\end{equation}
On the other hand, from the $u$-system \eqref{u_sym} with relations \eqref{u-v} and using condition \eqref{v_1},
\begin{align*}
& u_1' = 0 \Leftrightarrow \frac{v_1v_0}{1-v_1} + (2+r)v_0 +v_1 -\frac{r+2}{r+1}=0\\
& u_2' = 0 \Leftrightarrow (r+2)v_0 + v_1 - \frac{r+2}{r+1} + \frac{v_0v_1}{1-v_1}\\
& u_3'= 0 \Leftrightarrow  (1-v_1)^2 +(1-v_1)(\frac{1}{r+1}-(r+1)v_0) -v_0 = 0
\end{align*}
By solving the last line with respect to $(1-v_1)$, one has
$$1-v_1 = (r+1)v_0 \mbox{ or } 1-v_1 = -\frac{1}{r+1}.$$
But since $1-v_1 \in[0,1]$, necessarily $1-v_1 = (r+1)v_0$. Using \eqref{cond-v_1}, $v_0$ solves
$$2d(\lambda_1+\lambda_2r)(r+1)v_0^2 - (2d\lambda_1+2d\lambda_2r+r+1)v_0 + 1 = 0.$$
This implies
$$v_0 = \dfrac{1}{r+1} \mbox{ or } v_0 = \dfrac{1}{2d\lambda_1+2d\lambda_2r}.$$
\begin{enumerate}
\item if $v_0 = \frac{1}{r+1}$, then by relations \eqref{u-v}, 
$$u_1 = 0, \ u_3 = 1, \ v_2 = 1+u_2,$$
which is a contradiction.
\item if $v_0 = \dfrac{1}{2d\lambda_1+2d\lambda_2r}$, then
$$v_1 = \dfrac{r+1}{2d(\lambda_1+\lambda_2r)}, \ v_2 = \dfrac{r}{r+1}.$$
\end{enumerate}
Verifying this $v$-system to be a set of densities by condition \eqref{cond-01}, one case highlights a condition on $r$: $v_1 \leq 1$ if and only if $r(1-2d\lambda_2) \leq 2d\lambda_1-1$. In the case where $\lambda_2 \leq 1/(2d)$, then one has the condition 
\begin{equation}\label{cond-r_1}
r \leq \dfrac{2d\lambda_1-1}{1-2d\lambda_2}.
\end{equation}
Consequently, a non-trivial equilibrium of \eqref{v_sym} is given by 
\begin{equation}
(v_0,v_1,v_2) = \Big( \dfrac{1}{2d\lambda_1+2d\lambda_2r}, \dfrac{r+1}{2d\lambda_1+2d\lambda_2r}, \dfrac{r}{r+1}\Big)
\end{equation}
To put in a nutshell, this survey of equilibria associated to both mean-field models, in symmetric and asymmetric case, gave us the bounds \eqref{cond-r_0} and \eqref{cond-r_1} for the phase transition.

\section{Quenched random environment on $\mathbb Z^1$}\label{sec:quench}
We now derive bounds on the phase transition by using T.M. Liggett's results \cite{Lig91,Lig92} for the one-dimensional contact process in quenched environment, that we adapt to our different dynamics.
Let $p \in (0,1)$, define a random environment $\omega \in \{0,1\}^{\mathbb Z}$ where each site $x \in \mathbb Z$ is either free (0) with probability $1-p$ or slowed-down (1) with probability $p$, independently of any other site.

\bigskip
The contact process in random environment (CPRE) $(\xi_t)_{t \geq 0}$ that we consider has state space $\{0,1\} ^{\mathbb{Z}}$ with quenched environment $\omega$. Let $\lambda_1$ and $\lambda_2$ be growth parameters such that 
\begin{equation}\label{quench:birth}
\lambda_2 \leq  \lambda_c < \lambda_1,
\end{equation}
where $\lambda_c$ is the critical rate of the one-dimensional  contact process. 
Subsequently,
\begin{equation}\label{quench:rate}
p = r/(r+1)
\end{equation}
 for $r \in (0,\infty)$, is the probability a site is slowed down. 
Deaths occur at rate 1. If $\{\omega(x), x \in \mathbb Z\}$ is chosen deterministic (hence inhomogeneous), we will refer to it as the \textsl{contact process in inhomogeneous environment} (CPIE). Both CPRE and CPIE are still monotone according to Section \ref{sec:order}.

\bigskip
Denote by $\mathbb P_{\lambda_1, \lambda_2, r}^\omega$ the distribution of the CPRE with parameters $(\lambda_1, \lambda_2, r)$ and environment $\omega$. For fixed parameters $\lambda_1$ and $\lambda_2$, simplify by $\mathbb P_r^\omega$.  For any $A \subset \mathbb Z$, define $\Xi_t^A := \{ x \in \mathbb Z : \xi_t^A (x) =1 \}$, where $\xi_t^A$ denotes the process at time $t$ started from the initial configuration $\xi_0 = \mathbf{1}_A$. If $A=\{0\}$, simplify by $\Xi_t \equiv \Xi_t^{\{0\}}$. 
For almost-every realization of $\omega$, the CPRE is said to \textsl{survive} if
$$\mathbb P^\omega_r( \forall t \geq 0, \ \Xi_t \neq \emptyset) > 0$$ 
and \textsl{die out} otherwise. 

Consider subsequently two kinds of random environment: one depending on the vertices and one depending on the edges of the graph. The following results are adapted from \cite[Theorems 1 and 3]{Lig91} and \cite[Theorem 1.4]{Lig92}. 

\subsection{Random growth on vertices}\label{RE:vertex}
If $\lambda_v(k)$ is the growth rate from site $k \in \mathbb Z$ (See Figure \ref{fig:vertex}): a birth at site $k$ occurs at rate $\lambda_v(k-1)$ if $k-1$ is occupied plus at rate $\lambda_v(k+1)$ if $k+1$ is occupied, where
\begin{equation}\label{quenched:growth_v}
\lambda_v(k) = \lambda_1 (1-\omega(k))+\lambda_2 \omega(k)
\end{equation}
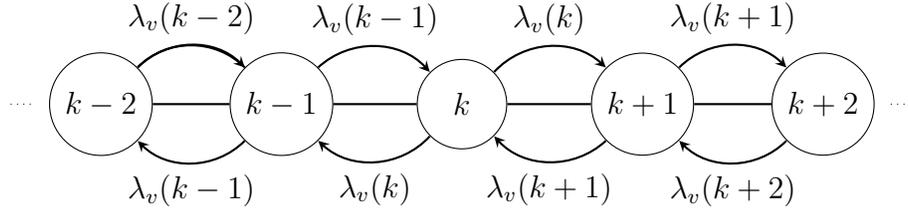
\begin{figure}[h]
\centering
\begin{tikzpicture} [domain=0:10, scale=0.6, baseline=0, >=stealth]
    \node (0) at (-4,0) [draw, circle={4pt}]{$k-2$};
    \node (1) at (0,0) [draw, circle={4pt}]{$k-1$};
	\node (2) at (4,0) [draw, circle={4pt}]{$ \ \ k \ \  $};
	\node (3) at (8,0) [draw, circle={4pt}]{$k+1$};
    \node (4) at (12,0) [draw, circle={4pt}]{$k+2$};
\draw[-,thick]    
	(0) edge [] node [below] {} (1)
	(1) edge [] node [below] {} (2)
	(2) edge [] node [below] {} (3)
	(3) edge [] node [below] {} (4)		
		(0) edge [bend left=45] node [sloped, anchor=left,above=0.2cm] {} (1);
\draw[dotted]	(13.5,0) -- (14,0);
\draw[dotted]	(-6,0) -- (-5.5,0);
\draw[->,thick]  
	(3) edge [bend left=45] node [above] {$\lambda_v(k+1)$} (4)
	(4) edge [bend left=45] node [below] {$\lambda_v(k+2)$} (3)
	(2) edge [bend left=45] node [above] {$\lambda_v(k)$} (3)
	(3) edge [bend left=45] node [below] {$\lambda_v(k+1)$} (2)
	(1) edge [bend left=45] node [above] {$\lambda_v(k-1)$} (2)
	(2) edge [bend left=45] node [below] {$\lambda_v(k)$} (1)
	(0) edge [bend left=45] node [above] {$\lambda_v(k-2)$} (1)	
	(1) edge [bend left=45, below] node [below] {$\lambda_v(k-1)$} (0);
	\end{tikzpicture}
	\caption{Random environment on vertices} \label{fig:vertex}
\end{figure}

\begin{theorem}\label{RE:vertex:extinct1}
Assume $\omega$ is a fixed environment. The CPIE dies out if for all $n \in \mathbb Z$,
\begin{equation}\label{inhomogeneous:extinction_eq}
\sum\limits_{k \geq n} \prod\limits_{j=n}^k \lambda_v (j+1) < \infty \mbox{ and } \sum\limits_{k \leq n} \prod\limits_{j=k}^n \lambda_v(j-1)  < \infty.
\end{equation}
\end{theorem}

If the family $\{\omega(k), k \in \mathbb Z\}$ is i.i.d. then the family $\{\lambda_v(k), k \in \mathbb Z\}$ is i.i.d as well and
\begin{corollary}\label{RE:vertex:extinct}
The CPRE dies out if $\mathbb E^\omega_r \log \lambda_v(0) < 0$.
That is, if $\lambda_2 < 1$ and $r > - \log \lambda_2 / \log \lambda_1$.
\end{corollary}
\begin{theorem}\label{RE:vertex:survival}
The CPRE survives if
$$\sum\limits_{j\geq 0} \mathbb E^\omega_r \Bigg( \dfrac{1}{\lambda_v(j)} \prod\limits_{k=1}^j \dfrac{\lambda_v(k) + \lambda_v(k-1)+ 1}{\lambda_v(k)\lambda_v(k-1)} \Bigg) < \infty.$$
\end{theorem}

No independence in the above product disables us to obtain explicit conditions for survival of the process. Nevertheless, by defining the randomness on the edges rather than on the vertices, meaning that the growth rates emanating from a site $k$ respectively to $k+1$ and to $k-1$ are randomly chosen for each $k \in \mathbb Z$, we are able to explicit bounds on $r$ with respect to $\lambda_1$ and $\lambda_2$.

\subsection{Random growth on oriented edges}\label{RE:edge}
\begin{figure}[h]
\centering
\begin{tikzpicture} [domain=0:10, scale=0.6, baseline=0, >=stealth]
    \node (0) at (-4,0) [draw, circle={4pt}]{$k-2$};
    \node (1) at (0,0) [draw, circle={4pt}]{$k-1$};
	\node (2) at (4,0) [draw, circle={4pt}]{$ \ \ k \ \  $};
	\node (3) at (8,0) [draw, circle={4pt}]{$k+1$};
    \node (4) at (12,0) [draw, circle={4pt}]{$k+2$};
\draw[-,thick]    
	(0) edge [] node [below] {} (1)
	(1) edge [] node [below] {} (2)
	(2) edge [] node [below] {} (3)
	(3) edge [] node [below] {} (4)		
		(0) edge [bend left=45] node [sloped, anchor=left,above=0.2cm] {} (1);
\draw[dotted]	(13.5,0) -- (14,0);
\draw[dotted]	(-6,0) -- (-5.5,0);
\draw[->,thick]  
	(3) edge [bend left=45] node [above] {$\lambda_e(k+2)$} (4)
	(4) edge [bend left=45] node [below] {$\rho_e(k+1)$} (3)
	(2) edge [bend left=45] node [above] {$\lambda_e(k+1)$} (3)
	(3) edge [bend left=45] node [below] {$\rho_e(k)$} (2)
	(1) edge [bend left=45] node [above] {$\lambda_e(k)$} (2)
	(2) edge [bend left=45] node [below] {$\rho_e(k-1)$} (1)
	(0) edge [bend left=45] node [above] {$\lambda_e(k-1)$} (1)	
	(1) edge [bend left=45, below] node [below] {$\rho_e(k-2)$} (0);
	\end{tikzpicture}
	\caption{Random environment on oriented edges} \label{fig:edges}
\end{figure}
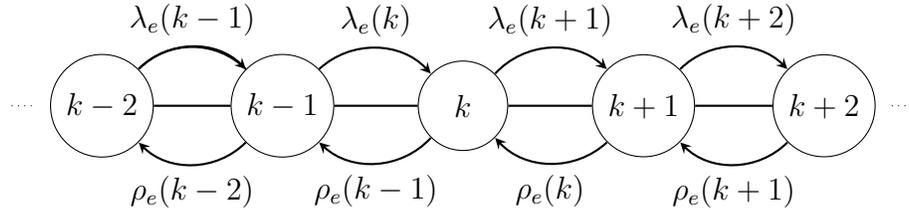

Let $\{(\rho_e(k), \lambda_e(k)), k \in \mathbb Z\}$ be an ergodic, stationary and i.i.d. sequence. Here, given a site $k \in \mathbb Z$, a birth from $k$ to $k+1$ occurs at rate  $\lambda_e(k+1)$ and independently of a birth from $k$ to $k-1$ occuring at rate $\rho_e(k-1)$. See Figure \ref{fig:edges}.

Suppose both rates are i.i.d, defined for all $k \in \mathbb Z$ by
\begin{equation}\label{quenched:growth_e}
\lambda_e(k+1) \stackrel{(d)}{=}\rho_e(k-1) \stackrel{(d)}{=}  \lambda_1 (1- \omega(k)) + \lambda_2 \omega(k).
\end{equation}

\begin{theorem}\label{quench:edge:extinct}
The CPRE dies out if
$$\mathbb E_r^\omega \lambda_e(k) < 1
\mbox{ and } 1-\mathbb E_r^\omega \Big( \dfrac{1}{\lambda_e(k)} \Big)< \mathbb E_r^\omega \Big(\dfrac{1}{\lambda_e(k)}\Big) \Big( 1 - \mathbb E_r^\omega \lambda_e(k) \Big). $$
that is, if $\lambda_2 < 1$ and $r > \dfrac{\lambda_1 - 1}{1 - \lambda_2}$.
\end{theorem}

\begin{theorem}\label{RE:edge:survival}
The CPRE survives if for all  $j \geq 0$,
$$\mathbb E^\omega_r \Big( \dfrac{1}{\lambda_e(j+1)} \Big) \Big( \mathbb E^\omega_r  \dfrac{\lambda_e(0) + \rho_e(0)+ 1}{\lambda_e(0)\rho_e(0)} \Big)^j < 1$$
that is, if $\lambda_2 < 1+\sqrt{2} < \lambda_1$ and 
 $r < {\lambda_2 \big(\lambda_1 - \sqrt{2} - 1 \big)}/\big({\lambda_1 \big( \lambda_2- \sqrt{2}-1\big) }\big)$.
\end{theorem}
\subsection{Numerical bounds on the phase transition}\label{RE:num}
Back to the basic contact process, assume $r =0$ then for all $x \in \mathbb Z$, $\omega(x) = 0$ a.s. and $\lambda_e(x)=\rho_e(x)=\lambda_1$. We thus recover the one-dimensional basic contact process with growth rate $\lambda_1$. In this case, our estimates lead to the following bound for $\lambda_c$.
\begin{corollary}\label{RE-cp}
For the one-dimensional basic contact process,
$$\lambda_c \leq 1+\sqrt{2}.$$
\end{corollary}
\begin{proof}
 According to Theorem \ref{RE:edge:survival}, if $r=0$ the process survives if 
$\lambda_1 > 1 + \sqrt{2}.$
\end{proof}
This bound is quite rough but we derived it simply. Consequently, one first deduces a bound on the critical value $\lambda_c$ on $\mathbb Z$: $\lambda_c \leq 1 + \sqrt{2} \simeq 2.41$. 

\bigskip
From previous results, let us derive numerical bounds on the segment where the phase transition occurs. Choosing $\lambda_1$ and $\lambda_2$ given by Theorem \ref{quench:edge:extinct} gives us lower bounds on the phase transition and satisfying condition from Theorem \ref{RE:edge:survival} gives us upper bounds.
\begin{center}
\begin{small}
$
\begin{array}{|c|c|c|}
\hline
\lambda_1 & \lambda_2 & \mbox{phase transition}\\
   \hline
1000&0.8& [0.49,4995)\\
100&0.8&[0.48,495)\\
10&0.8&[0.36,45]\\
2 &0.8& (0,5] \\
\hline
\end{array}
$
\hspace*{1cm}
$
\begin{array}{|c|c|c|}
\hline
\lambda_1 & \lambda_2 & \mbox{phase transition}\\
   \hline
1000&1.4& [1.37,\infty)\\
100&1.4&[1.34,\infty)\\
10&1.4&[1.04,\infty)\\
2 &1.4& \mathbb R_+ \\
\hline
\end{array}
$
\end{small}
\end{center}

Remark that the necessary condition $\lambda_2 < 1$ disables us to conclude to an upper bound for values of $\lambda_2$. In a similar way, Theorem \ref{RE:edge:survival} imposes $\lambda_1$ to be larger than $1+\sqrt{2}$,  disabling us to find an explicit lower bound in such cases.

\subsection{Proofs}
The \textsl{contact process in random environment} $(\xi_t)_{t \geq 0}$ on $\{0,1\}^\mathbb Z$ introduced by T.M. Liggett \cite{Lig91,Lig92} has transitions at $x \in \mathbb Z$ given by
\begin{equation}
\begin{array}{ll}
& 0 \rightarrow 1 \mbox{ at rate }  \rho(x)\xi(x+1) + \lambda(x) \xi(x-1)\\ 
& 1 \rightarrow 0 \mbox{ at rate } \delta(x)
\end{array}
\end{equation}
where the family $\{\big( \delta(x),\rho(x),\lambda(x)\big), x \in \mathbb Z\}$ stands for the random environment, which is an ergodic stationary process. To apply \cite[Theorems 1 and 3]{Lig91} and \cite[Theorem 1.4]{Lig92} to our case, the link to \eqref{quenched:growth_v} and \eqref{quenched:growth_e} is given by: for all $k \in \mathbb Z$,
$$\lambda(k+1) = \rho(k-1) = \lambda_v(k)  = \rho_v(k)$$
$$\lambda(k) = \lambda_e(k) \mbox{ and } \rho(k) = \rho_e(k)$$
\begin{proof}[Proof of Corollary \ref{RE:vertex:extinct}]
By the ergodic theorem, $$\lim\limits_{k \rightarrow \infty} \dfrac{1}{k} \sum\limits_{j=0}^k \log \lambda_v(j) = \mathbb E_r^\omega \log_v(0),$$
where $E_r^\omega \log_v(0) = \dfrac{r\log \lambda_2 + \log \lambda_1}{r+1}$, so that  by Theorem \ref{RE:vertex:extinct1} extinction occurs if $$r > - \log \lambda_1 / \log \lambda_2.$$
\end{proof}

\begin{proof}[Proof of Theorem \ref{quench:edge:extinct}]
Since $\mathbb E_r^\omega \lambda_e(k) = \lambda_1(1-p)$, we have both conditions
\begin{equation}\label{condRE}
\lambda_2 < 1 \mbox{ and } r > (\lambda_1-1)/ (1-\lambda_2).
\end{equation}
Meanwhile, the second assertion of Theorem \ref{quench:edge:extinct} holds if the polynomial 
$$A(r):= 2 r^2 \dfrac{1-\lambda_2}{\lambda_2} + r \Big( \dfrac{2-\lambda_2}{\lambda_1} + \dfrac{2-\lambda_1}{\lambda_2} - 2 \Big) + 2\dfrac{1-\lambda_1}{\lambda_1}$$
is positive which holds if 
$$r > \dfrac{(\lambda_1+\lambda_2-2)(\lambda_1+\lambda_2) +  (\lambda_1-\lambda_2)\sqrt{(\lambda_1+\lambda_2-2)^2+4\lambda_1\lambda_2}}{4\lambda_1(1-\lambda_2)},$$
but this condition is satisfied as soon as \eqref{condRE} are.
\end{proof}

\begin{proof}[Proof of Theorem \ref{RE:edge:survival}]
Condition of Theorem \ref{RE:edge:survival} is satisfied if the following polynomial 
\begin{multline*}
A(r) := r^2 \Big[  \lambda_1^2(2\lambda_2+1)-\lambda_1^2 \lambda_2^2 \Big] + r \lambda_1 \lambda_2 \Big[ 2(\lambda_1 + \lambda_2 + 1) - 2 \lambda_1 \lambda_2 \Big]  + \Big[ \lambda_2^2 (2\lambda_1 +1) - \lambda_1^2 \lambda_2^2  \Big] 
\end{multline*}
is negative, which holds if $$\lambda_1 > 1 +\sqrt 2, \ \lambda_2 < 1 + \sqrt 2 \mbox{ and } r < \dfrac{- \lambda_2 \Big(\lambda_1 - \sqrt{2} - 1 \Big)}{\lambda_1 \Big( \lambda_2- \sqrt{2}-1\Big) }.$$
Since we assumed $\lambda_2 < \lambda_c$, where $\lambda_c \leq 1 + \sqrt 2$ by Corollary \ref{RE-cp}, only conditions of Theorem \ref{RE:edge:survival} remain.
\end{proof}

\subsection*{Acknowledgements} \textit{I would like to thank Rinaldo Schinazi and Ellen Saada for having introduced me to this topic, motivating discussions and uncountable comments on previous versions.}


\end{document}